\documentclass{elsarticle}

\usepackage{amsmath}
\usepackage{amssymb}
\usepackage{bm}
\usepackage{tikz}
\usepackage{pgfplots}

\usepackage{hyperref}

\newproof{pf}{Proof}

\journal{arXiv} 

\graphicspath{{figs/}}

\begin{document}

\begin{frontmatter}

\title{Numerical solution of BVP for the incompressible Navier-Stokes equations at large Reynolds numbers\tnoteref{label1}}
\tnotetext[label1]{The work was supported by the Russian Science Foundation (grant No. 23-41-00037).}

\author[msu]{D.V.~Lomasov}
\ead{lomasovdv@my.msu.ru}

\author[msu,svfu]{P.N.~Vabishchevich\corref{cor}}
\ead{vab@cs.msu.ru}

\address[msu]{Lomonosov Moscow State University, 1, Leninskie Gory,  Moscow, Russia}
\address[svfu]{North-Eastern Federal University, 58, Belinskogo st, Yakutsk, Russia}

\cortext[cor]{Corresponding author}

\begin{abstract}

The problems of numerical modeling of viscous incompressible fluid flows are widely considered in computational fluid dynamics.
Stationary solutions of boundary value problems for the Navier-Stokes equations exist at large Reynolds numbers, but they are unstable and lead to transient or turbulent unsteady regimes.
In addition, the solution of the boundary value problem at large values of Reynolds number may be non-unique.
In this paper, we consider computational algorithms numerical algorithms for finding such stationary solutions.
We use natural pressure-velocity variables under standard finite element approximation on triangular grids.
Iterative methods with different linearizations of convective transport are used to test a two-dimensional problem of incompressible fluid flow in a square-section cavity with a movable top lid.
The developed computational algorithm allowed us to obtain two solutions when the Reynolds number exceeds a critical value
for flows in a cavity of semi-elliptical cross-section.

\end{abstract}


\end{frontmatter}

\section{Introduction}

Applied continuum mechanics often uses the Navier-Stokes equations to model viscous incompressible fluids, essential in fluid dynamics \cite{Batchelor,LandauLifshic1986}.
The flow is typically laminar at low Reynolds numbers, exhibiting smooth and stable motion.
As the Reynolds number increases beyond a critical point, the flow becomes unsteady, indicating the onset of time-dependent behavior.
With further increases in the Reynolds number, the flow transitions to a turbulent state characterized by chaotic and irregular motion.
While the Navier-Stokes equations are foundational, they have limitations in fully capturing the complexities of turbulence, especially at smaller scales where thermal fluctuations are significant.

In fluid dynamics, stationary solutions of the Navier-Stokes equations are present for all Reynolds numbers \cite{lions1}.
As the Reynolds number increases, particularly at elevated values, the uniqueness of these stationary solutions may become an issue.
In regimes characterized by high Reynolds numbers, which are commonly linked to turbulent flows, multiple stationary solutions can arise from the same boundary conditions.

A promising area in computational fluid dynamics is the development of algorithms to identify and analyze multiple stationary solutions for viscous incompressible fluid flows at high Reynolds numbers.
This challenge arises from potential nonlinearities and instabilities in the governing equations.
Addressing these complexities requires consideration of the singularities in the boundary value problem, which can occur due to rapid flow changes or sharp gradients near boundaries.
To achieve accurate solutions, we must use sufficiently detailed computational grids.
Finer grids enhance flow feature resolution and numerical stability, capturing the nuances of fluid behavior and improving the reliability of computed stationary solutions.
Exploring stationary solutions at large Reynolds numbers emphasizes the intricate nature of fluid dynamics and the necessity for advanced computational methods to uncover multiple solutions.

When tackling boundary value problems for partial differential equations in irregular domains, finite element, and finite volume methods are the most widely used numerical techniques \cite{KnabnerAngermann2003,QuarteroniValli}.
Both methods effectively model incompressible flows and can accommodate complex geometries and varying boundary conditions \cite{glowinski2022numerical,roychowdhury}.
For two-dimensional flow problems, using vorticity and stream function variables instead of the standard pressure-velocity formulation can simplify the governing equations and enhance numerical stability. This method effectively captures flow characteristics while reducing the complexity of pressure calculations.

Computational algorithms are crucial in addressing challenges in fluid dynamics and are typically validated through benchmark test cases.
A widely recognized benchmark for the Navier-Stokes equations related to incompressible flows is the lid-driven cavity problem \cite{cavity,kuhlmann2019lid}.
This scenario consists of a square cavity with a top lid that moves constantly, generating an extensively studied flow pattern. It is a reference for evaluating numerical methods \cite{roychowdhury}.
The lid-driven cavity problem is significant due to its complex flow dynamics, especially at elevated Reynolds numbers \cite{erturk2005numerical,hostos2021solving,wahba2012steady}.
Recent studies have successfully determined numerical solutions for two-dimensional steady-state flow within a square cavity at Reynolds numbers of  $10\,000$ and higher.
These high Reynolds number cases are essential for comprehending turbulent flow behaviors commonly observed in real-world applications.
We can use various iterative techniques to derive approximate solutions for these challenges.
Such methods facilitate the refinement of solutions through successive approximations.
Nonetheless, it is crucial to recognize that employing basic methods for stationary problems may yield ineffective results.
When the Reynolds number surpasses a critical value, the flow transitions to a nonstationary state, potentially exhibiting periodic behavior, complicating the solution process.
Recent developments in the area have identified two stationary solutions for a semi-elliptical cavity \cite{erturk2022bifurcation}.
This discovery underscores the complexity and diversity of flow behaviors in various geometrical configurations, further highlighting the necessity for robust computational methods to investigate these phenomena.
Numerical results for these fluid dynamics problems are often obtained using the vorticity-stream function formulation.

This study focuses on solving two-dimensional fluid dynamics problems using natural variables, specifically pressure, and velocity, along with standard finite element approximations on triangular grids.
This approach is well-suited for capturing the complexities of fluid flow, particularly in geometrically intricate domains.
We introduce a relaxation-type computational algorithm based on a specialized linearization of the convective term in the Navier-Stokes equations.
This linearization is crucial for managing the nonlinearity associated with convective transport in fluid flows, particularly at high Reynolds numbers where flow behavior can become chaotic and complex.
The relaxation algorithm iteratively refines the solution, allowing convergence towards an approximate solution regardless of the initial guess.
This characteristic is particularly beneficial when dealing with flows characterized by large Reynolds numbers, where traditional methods may need help finding a stable solution.
Our approach explicitly applies this algorithm to the flow problem in a driven semi-elliptical cavity.
This scenario presents unique challenges due to the geometry and the potential for complex flow patterns.
By carefully selecting the relaxation parameter within our algorithm, we can derive two distinct solutions for the flow within the cavity.
This ability to obtain multiple solutions highlights the flow dynamics' sensitivity to initial conditions and parameter choices, which is a critical aspect of fluid dynamics research.

\section{Problem statement} 

Our study addresses the problem of finding a numerical solution for the steady-state problem for the equations of viscous incompressible fluid. In the bounded domain $\Omega$, the dimensionless Navier-Stokes equations are used:
\begin{equation}\label{1}
C(\bm u, \bm u) + \operatorname{grad} p - \frac{1}{\operatorname{Re}} \operatorname{div} \operatorname{grad} \bm u = 0 ,
\end{equation}
\begin{equation}\label{2}
\operatorname{\operatorname{div}}  \bm u=0,
\quad \bm x \in \Omega .
\end{equation}
Here, $\bm u$ is the fluid velocity, $\operatorname{Re}$ is the Reynolds number, and $p$ is the pressure normalized to the density.
At the boundary, the velocity is considered to be prescribed:
\begin{equation}\label{3}
\bm u(\bm x) = \bm g(\bm x),
\quad \bm x \in \partial \Omega .
\end{equation}
The convective transport in equation (\ref{1}) is taken in conservative form:
\begin{equation}\label{4}
C(\bm v, \bm u) = \operatorname{div}  (\bm v\otimes\bm u) .
\end{equation}

The boundary value problem is nonlinear because there is a convective term $C(\bm u, \bm u)$ in the equation.
At small Reynolds numbers, the nonlinearity of the boundary value problem (\ref{1})--(\ref{4}) is weakly manifested, and therefore, there are no problems in finding an approximate solution.
We are interested in the case of large values of the Reynolds number.
In this case, we must take into account the singularity of the boundary value problem under consideration, which is generated by the small parameter $1/\operatorname{Re}$ at higher derivatives in equation (\ref{1}).
Due to this circumstance, we must use sufficiently detailed computational grids.
At large $\operatorname{Re}$ values, the problem of approximate solution of the boundary value problem (\ref{1})--(\ref{3}) is enormously complicated by the nonlinearity of the equation. Therefore, the critical point in developing the computational algorithm is related to the choice of the iterative method.

At large Reynolds numbers, the solution of the stationary problem for the incompressible fluid equations may be non-unique.
The problem of isolation of separate solutions requires even more attention.
At present, we have only the first results in this direction.
We want robust computational algorithms for obtaining multiple solutions of stationary Navier-Stokes equations.

The object of our research is the construction of iterative methods.
At each iteration, we use well-established computational fluid dynamics methods.
Considering the necessity of solving problems in irregular computational domains, we apply standard finite element approximations.
In computational practice on triangular (tetrahedral) grids, the Hood-Taylor finite element \cite{glowinski2022numerical} for the Navier-Stokes equations for viscous incompressible fluid is most often used.
In this case, the velocity is approximated using second-degree polynomials, and the pressure --- using first-degree polynomials.
The research software was written in Python using the extensive capabilities of the \textsf{FEniCS} \cite{logg2012automated} computing platform.
All calculations were performed on personal computing equipment.
In our two-dimensional calculations, we used a direct solver to solve the system of linear equations.
Computational experiments were performed with support for parallel computing.
The program \textsf{gmsh} \cite{geuzaine2009gmsh} is used to construct the triangular computational grid.

\section{Newton's method} 

Newton's method is used in the numerical solution of nonlinear problems.
The possibilities of such a standard approach are illustrated by the results of numerical solution of a two-dimensional problem of incompressible fluid flow in a cavity of square cross-section.
The calculations are performed on the basis of gradual increase of Reynolds number and using the obtained solution as an initial one to obtain an approximate solution at large Reynolds numbers.

\subsection{Iterative method} 

Let the approximate solution of $\bm u^k$ at $k$ iterations be known, the initial approximation $\bm u^0$ is given.
When using the two-level iterative method, the new approximation is obtained from the solution of the problem
\begin{equation}\label{5}
 \big (C (\bm u, \bm u) \big )^{k+1} + \operatorname{grad} p^{k+1} - \frac{1}{\operatorname{Re}} \operatorname{div}
 \operatorname{grad} \bm u^{k+1} = 0 ,
\end{equation}
\begin{equation}\label{6}
 \operatorname{\operatorname{div}}  \bm u^{k+1} =0,
 \quad \bm x \in \Omega ,
\end{equation}
\begin{equation}\label{7}
 \bm u^{k+1}(\bm x) = \bm g(\bm x),
 \quad \bm x \in \partial \Omega .
\end{equation}
The type of iterative method is related to the choice of approximations of the convective term.
Taking into account the quadratic nonlinearity $C (\bm u, \bm u)$, for Newton's method we have
\begin{equation}\label{8}
 \big (C (\bm u, \bm u) \big )^{k+1} = C (\bm u^{k}, \bm u^{k+1}) + C (\bm u^{k+1}, \bm u^{k}) - C (\bm u^{k}, \bm u^{k}) .
\end{equation} 

The convergence of the iterative method was evaluated by the difference of the approximate solution between iteration steps 
\[
\varepsilon (k) = \|\bm u^{k}(\bm x) - \bm u^{k-1}(\bm x)\| ,
\quad k = 1,2, \ldots .
\]
Here $\| \cdot\|$ is the norm in $\bm L_2(\operatorname{\Omega})$:
\[
 \|\bm u (\bm x) \| = \left ( \int_{\Omega} \bm u^2(\bm x) \, d \bm x \right)^{1/2} . 
\] 

\subsection{Square lid-driven cavity flow}

The above iterative methods were used to approximate the problem of viscous incompressible fluid flow in a cavity. 
The stationary boundary value problem for the incompressible Navier-Stokes equations (\ref{1})--(\ref{4}) is solved inside a square cavity
\[
\Omega  = \{ \bm x \ | \ \bm x = (x_1, x_2), \ 0 < x_\alpha < 1, \ \alpha =1,2 \}.
\]
The boundary conditions (\ref{3}) have the form
\[
\bm g(\bm x) =
\left \{ \begin{array}{cc}
(1, 0),   & \ x_2 = 1, \\
(0, 0),     & \ x_2 < 1, \\
\end{array}
\right .
\quad \bm x \in \partial \Omega .
\]

The stream function is considered the main characteristic of the flow.
After solving the problem (\ref{1})--(\ref{4}) the stream function is calculated from the solution of the boundary value problem
\[
\operatorname{div} \operatorname{grad} \psi = -\frac{\partial u_1}{\partial x_2} +\frac{\partial u_2}{\partial x_1},
\quad \bm x \in \Omega ,
\]
\[
\psi (\bm x) = 0, \quad \bm x \in \partial \Omega .
\]
From the known velocity field, we also calculate the vorticity $w \ (\operatorname{rot} \bm u = (0,0, w))$.

\begin{figure}[htbp]
\begin{tabular}{cc}
    $m=64$  	& $m=128$  \\ 
    \includegraphics[width=0.45\textwidth]{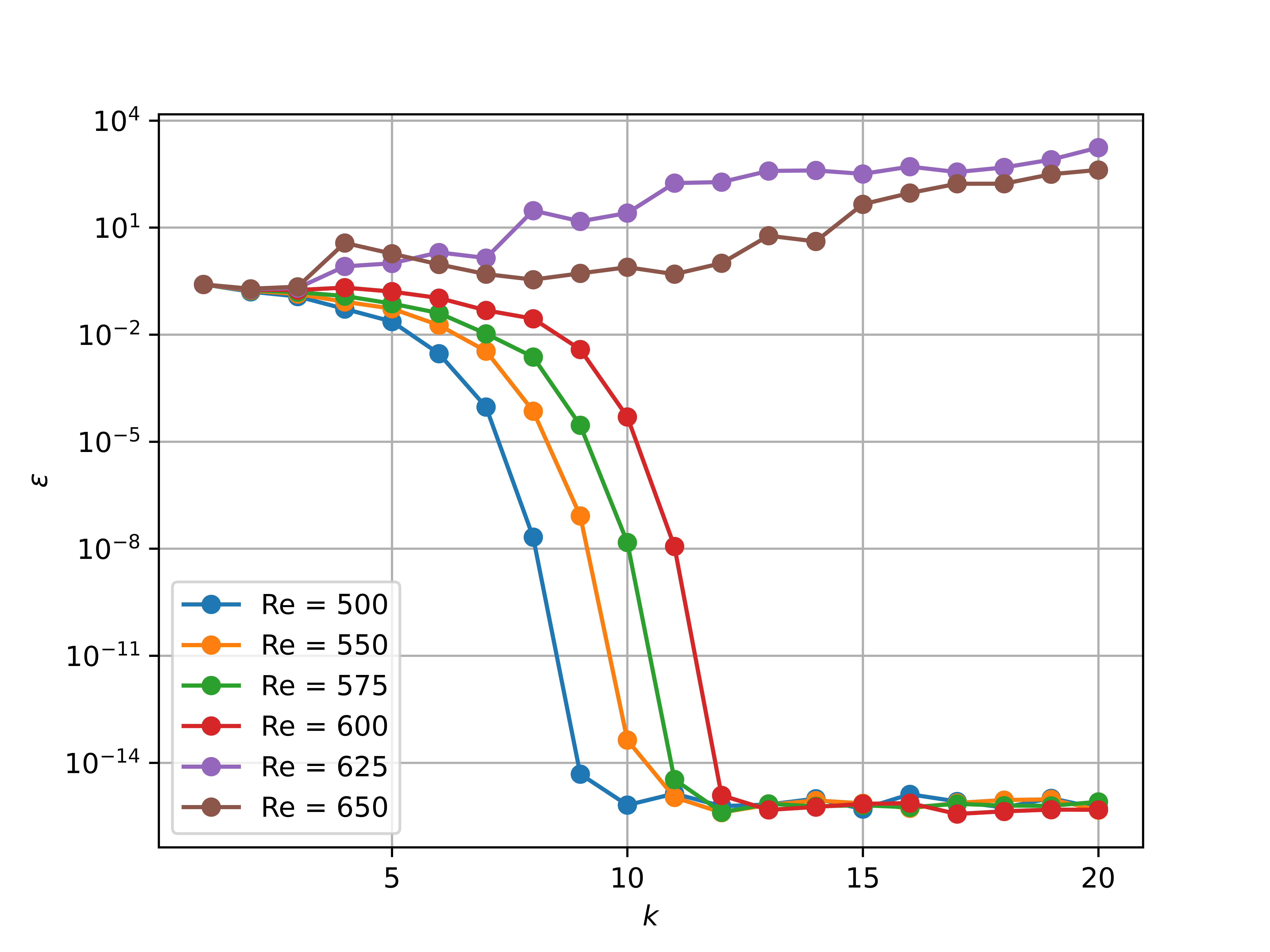} &
    \includegraphics[width=0.45\textwidth]{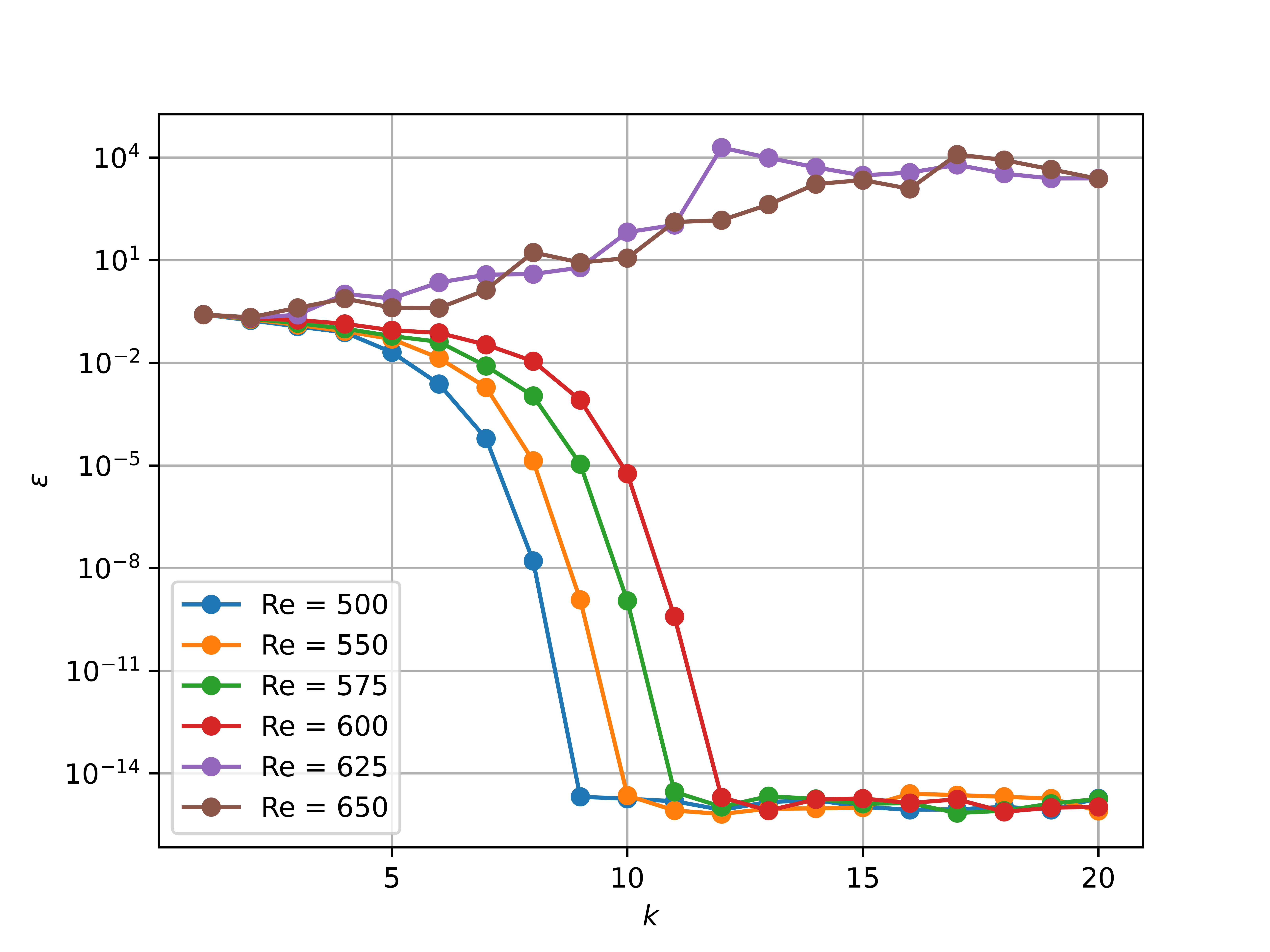}\\
    $m=256$  	&  $m=512$ \\ 
    \includegraphics[width=0.45\textwidth]{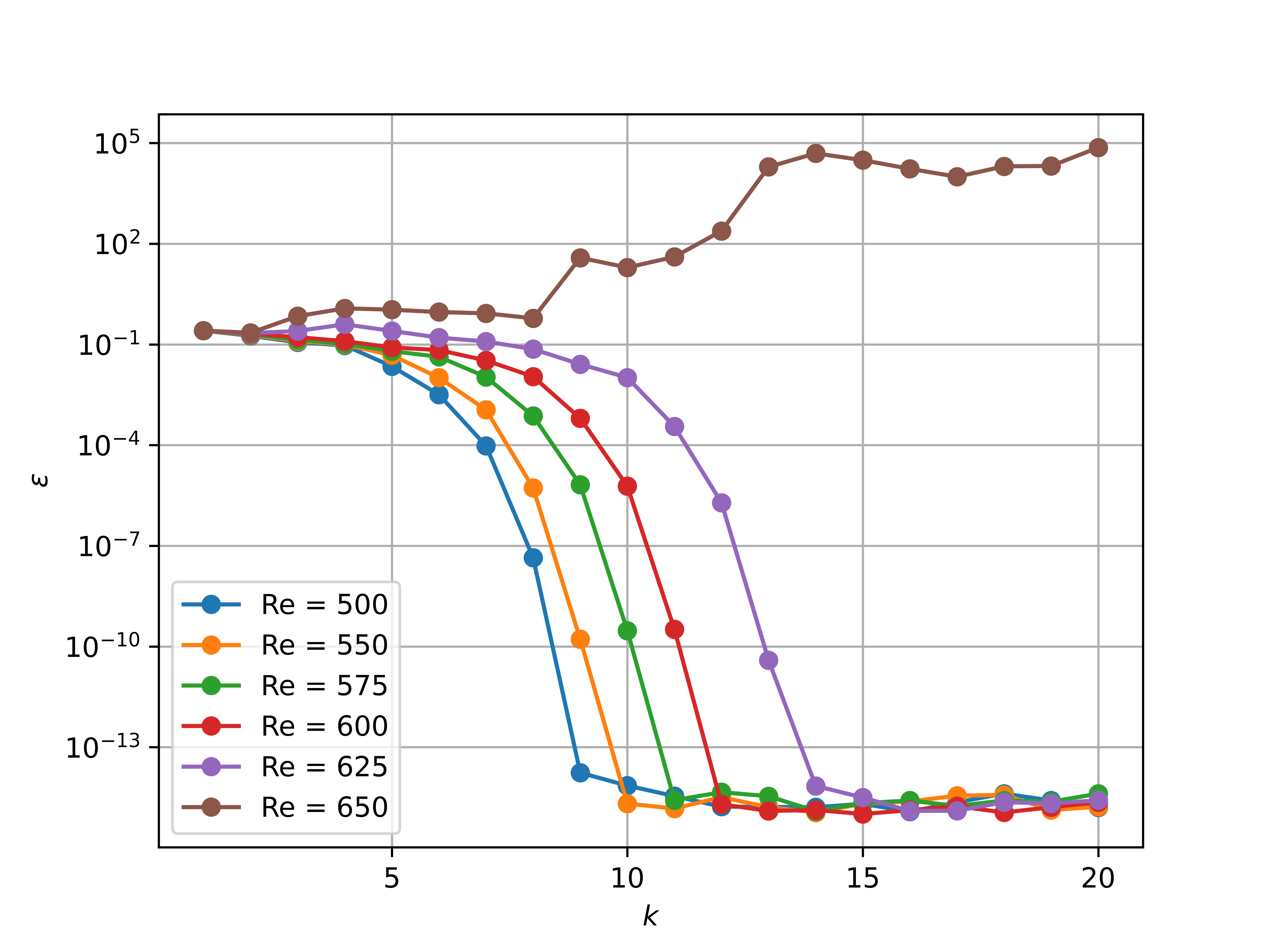} &
    \includegraphics[width=0.45\textwidth]{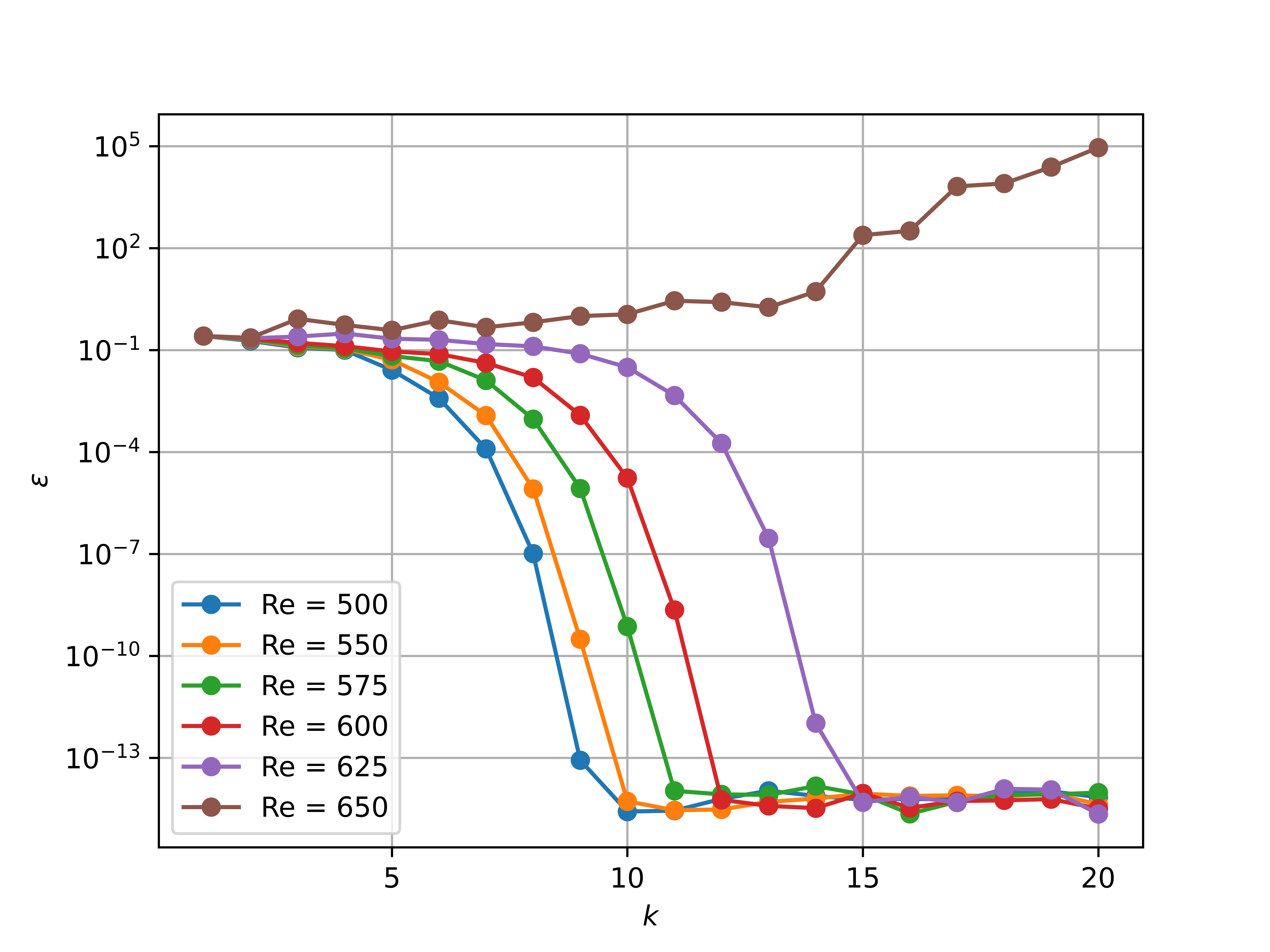}  
\end{tabular}
\caption{Convergence of Newton's method at different Reynolds numbers.}
\label{f-1}
\end{figure}

The flow in the cavity is modeled on a triangular grid with a uniform partition into $m$ intervals along each direction.
Fig.~\ref{f-1} shows the computational results when applying the iterative Newton method (\ref{5})--(\ref{8}) with the initial approximation $\bm u^0 = 0, \ \bm x \in \Omega$. 
The convergence of Newton's method is observed at small Reynolds numbers, at which the approximate solution is very fast.
The dependence on the used computational grid is weakly expressed.

\subsection{Other forms of convective transport}

In computational fluid dynamics, the problem of recording convective transport is actively discussed, particularly on the conservative (divergent) or characteristic (non-divergent) form of the equations of motion of a continuous medium.
When incompressibility is taken into account (equation (\ref{2})), the convective transport operator $C (\bm v, \bm u)$ instead of (\ref{4}) is often written in characteristic form
\begin{equation}\label{9}
C(\bm v, \bm u) = (\bm v \cdot \operatorname{grad}) \bm u = \operatorname{div} (\bm v\otimes\bm u) - \bm v \operatorname{div} \bm u.
\end{equation}
When applying the finite element method, the symmetric form is of particular note
\begin{equation}\label{10}
C(\bm v, \bm u) = \operatorname{div} (\bm v\otimes\bm u) - \frac{1}{2}  \bm v \operatorname{div} \bm u.
\end{equation}
In this case, we take the half-sum of the convective transfer operator in conservative form (\ref{4}) and characteristic form (\ref{9}).
The operator (\ref{10}) is skew-symmetric (energy neutral) in $\bm L_2(\operatorname{\Omega})$ for any $\bm v(\bm x)$.
Due to this fact, this choice of the convective transfer operator is said to be a skew-symmetric form.

\begin{figure}[htbp]
\begin{tabular}{cc}
$\operatorname{div}  (\bm u\otimes\bm u) - \bm u \operatorname{div} \bm u$  & $\operatorname{div}  (\bm u \otimes\bm u) - \frac{1}{2}  \bm u \operatorname{div} \bm u$  \\
\includegraphics[width=0.45\textwidth]{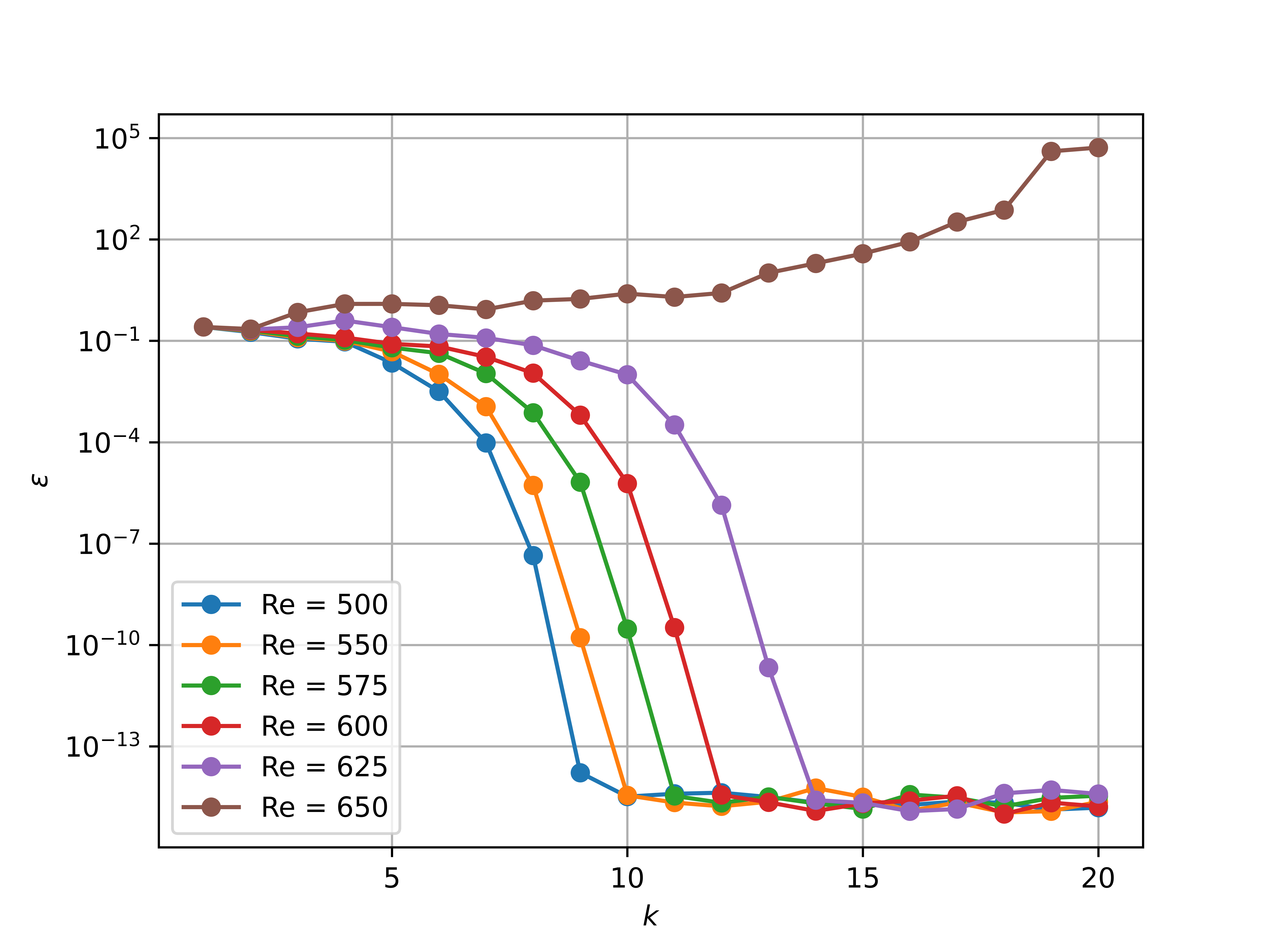} &
\includegraphics[width=0.45\textwidth]{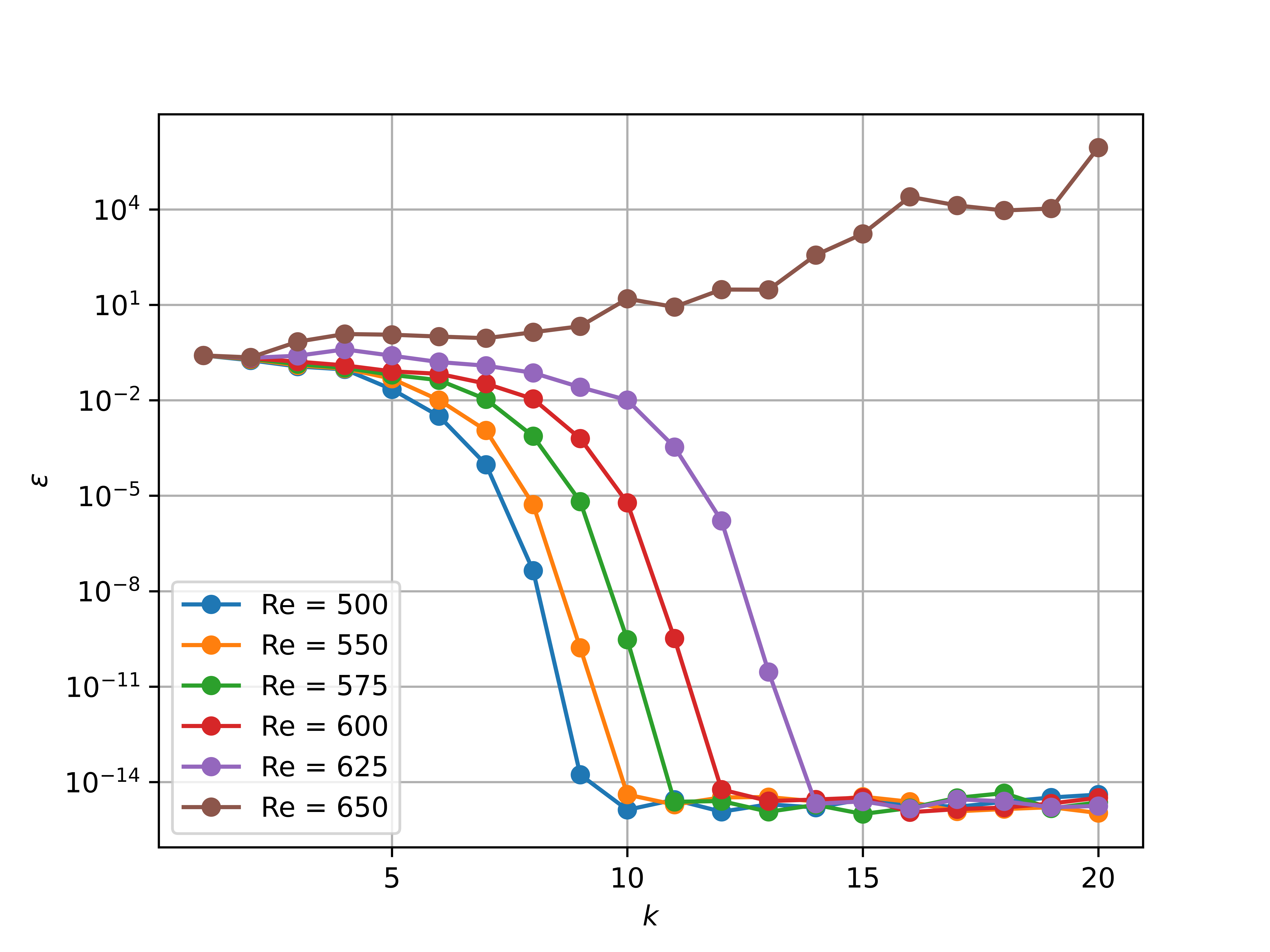}
\end{tabular}
\caption{Convergence of Newton's method under different forms of convective transport.}
\label{f-2}
\end{figure}

The convergence of Newton's method (\ref{5})--(\ref{8}) when the convective summand is given according to (\ref{9}) or (\ref{10}) is illustrated in Fig.~\ref{f-2}.
The calculations are performed on a uniform grid with $m=256$, and for the conservative form (\ref{4}), the data are shown in Fig.~\ref{f-1}.
No significant influence of the choice of form for convective transport $C (\bm v, \bm u)$ is observed.
The modifications of the iterative methods discussed below are performed by setting $C (\bm v, \bm u)$ according to (\ref{4}).

\subsection{Solution at large Reynolds numbers}

Newton's method for the approximate solution of nonlinear systems of equations converges, provided that the initial approximation is sufficiently close to the exact solution.
In the approximate solution of the boundary value problem (\ref{1})--(\ref{4}), a good approximation can be a solution at a smaller Reynolds number.
A natural strategy for obtaining a solution at large $\operatorname{Re}$ is to solve problems with sequentially
\[
\operatorname{Re}_1 < \operatorname{Re}_2 < \cdots < \operatorname{Re}_n < \cdots < \operatorname{Re}_N = \operatorname{Re}.
\]
The solution found at $\operatorname{Re}_{n-1}$ is taken as an initial approximation in the approximate solution of the problem with $\operatorname{Re}_{n}$.

\begin{figure}[htbp]
\begin{tabular}{cc}
$m=128$  & $m=256$  \\
\includegraphics[width=0.45\textwidth]{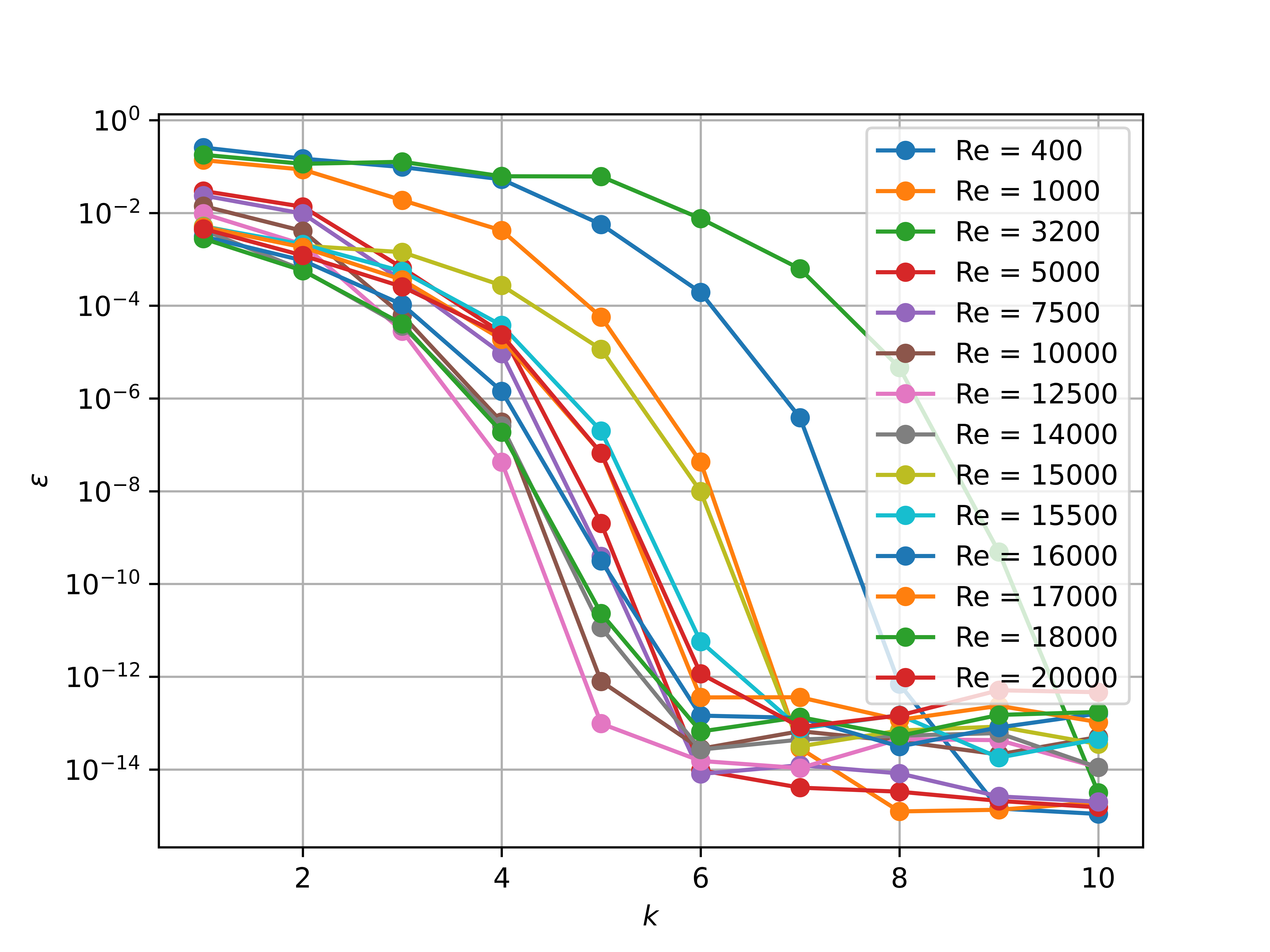} &
\includegraphics[width=0.45\textwidth]{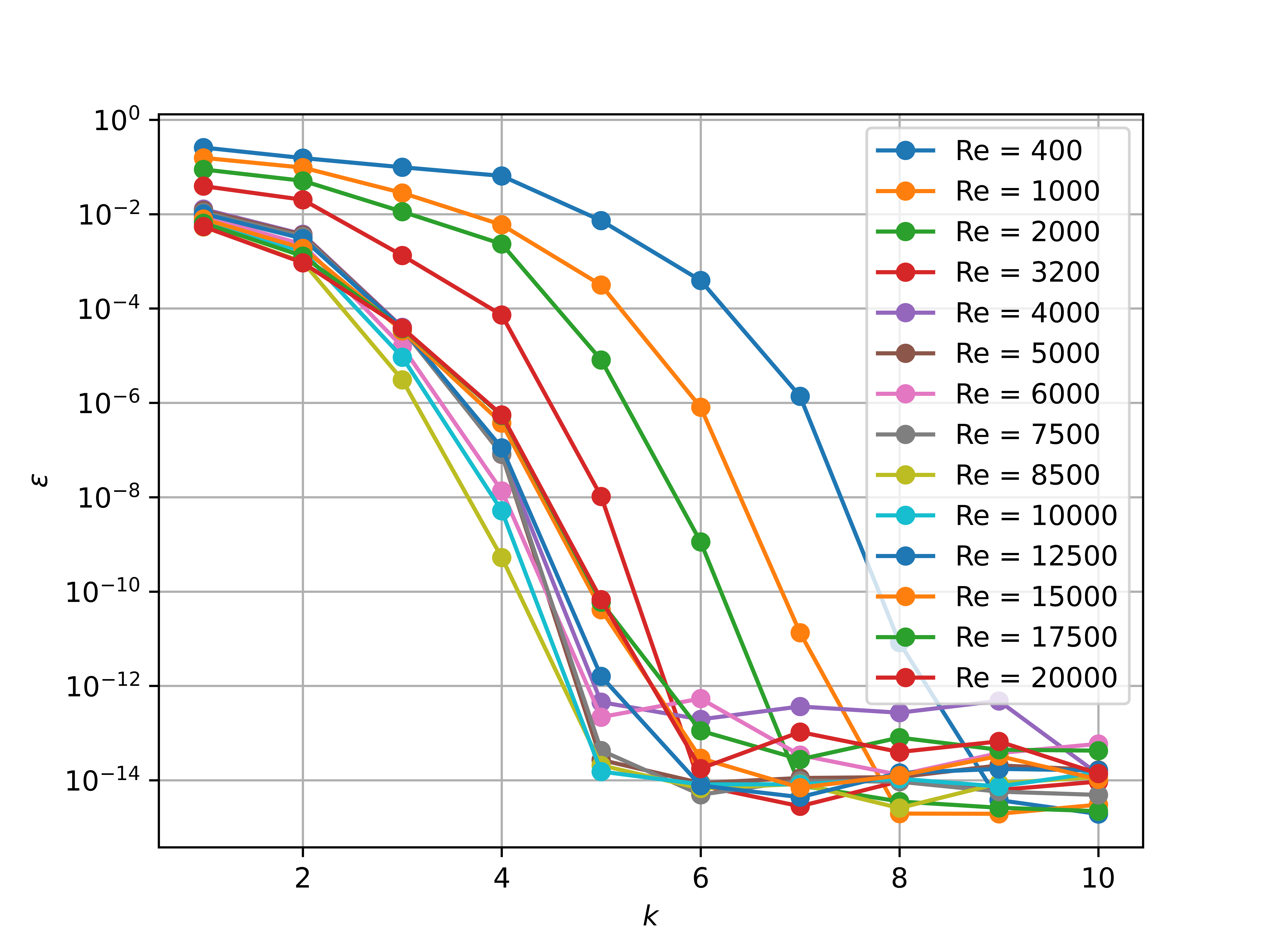}
\end{tabular}
\caption{A sequential increase in Reynolds number.}
\label{f-3}
\end{figure}

Fig.~\ref{f-3} shows the convergence of Newton's method for sequential increase of $\operatorname{Re}_{n}$.
An approximate solution is obtained when $\operatorname{Re} = 20\,000$ on two computational grids: $m = 128$ and $m = 256$.
At each $\operatorname{Re}_{n}$, 10 iterations were made in Newton's method.
In our case, the total number of iterations to obtain an approximate solution is $140$.

\begin{figure}[htbp]
\begin{tabular}{cc}
$m=64$  & $m=128$  \\
\includegraphics[width=0.45\textwidth]{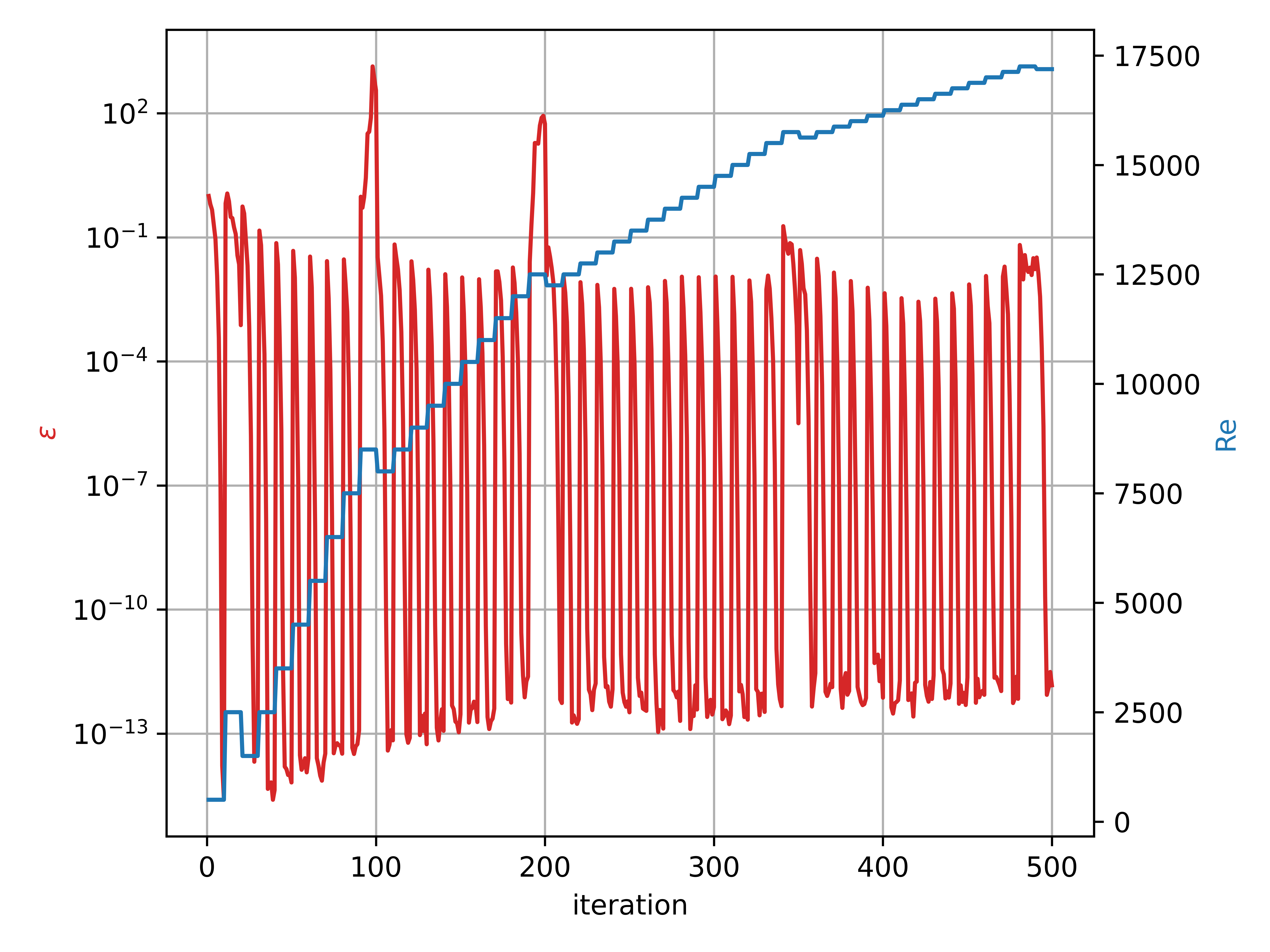} &
\includegraphics[width=0.45\textwidth]{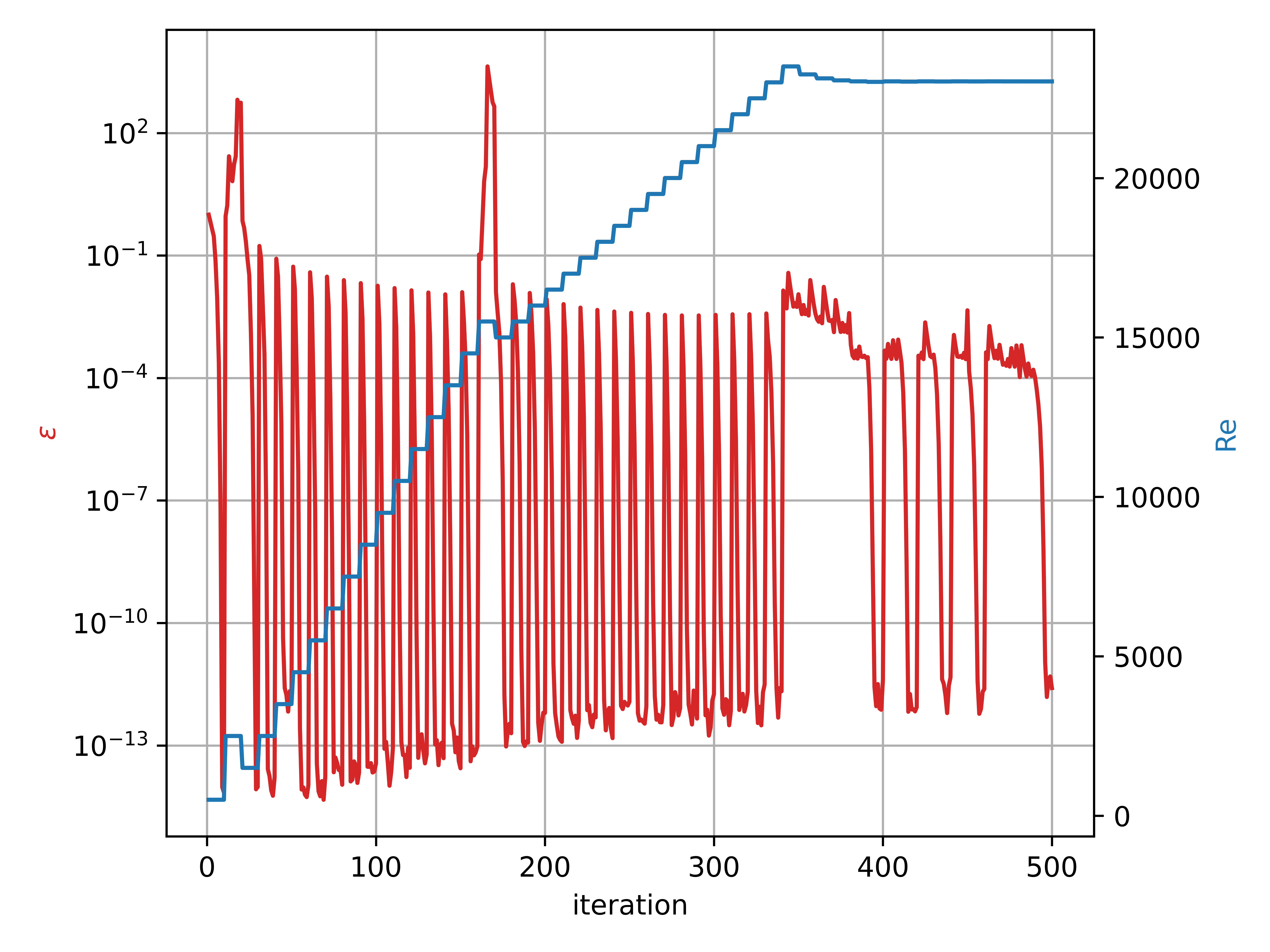} \\
$m=256$  &  $m=512$ \\
\includegraphics[width=0.45\textwidth]{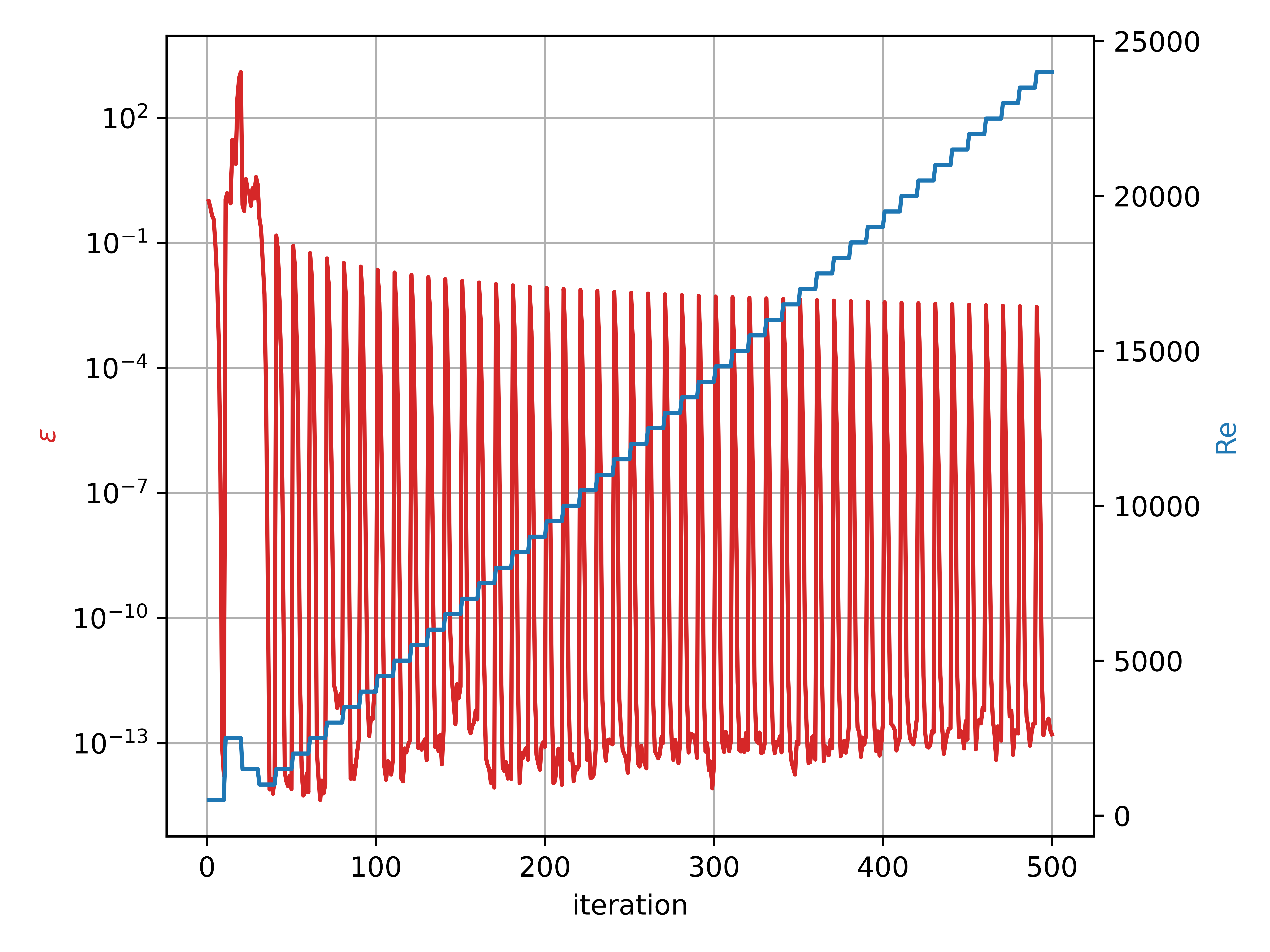} &
\includegraphics[width=0.45\textwidth]{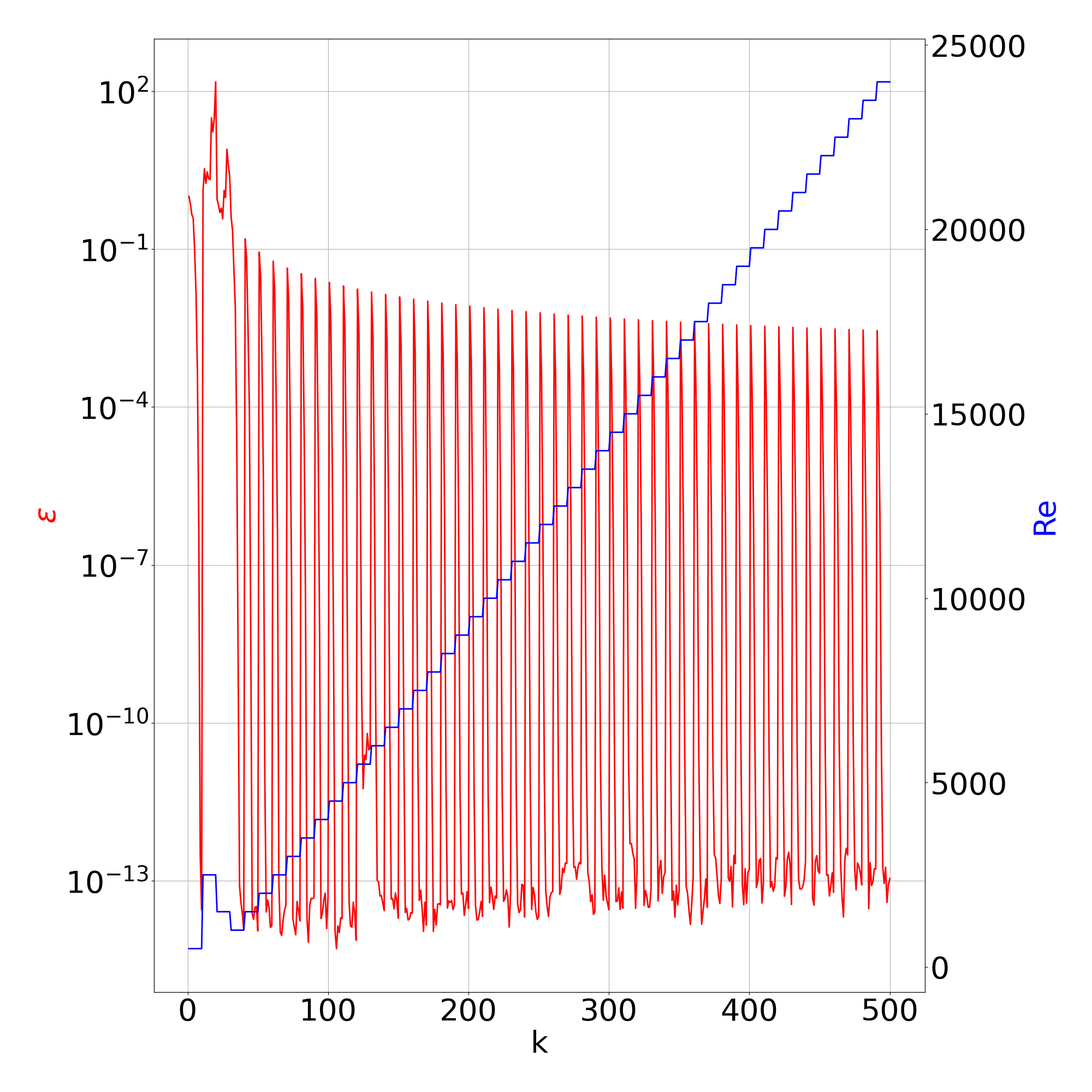}
\end{tabular}
\caption{Refinement of Reynolds number by bisection method.}
\label{f-4}
\end{figure}

In the sequential Reynolds number increase under consideration, the key is to select the sequence $\operatorname{Re}_{n}, \ n = 1,2, \ldots, N$.
We can automate this process by using the following bisection algorithm.
We will assume that there is a value of Reynolds number $\operatorname{Re}_{n}$ at which we have an approximate solution.
The increment (current increment step) of the Reynolds number $\delta_n$ is given.
The new trial value is defined according to $\operatorname{Re}_{n+1} = \operatorname{Re}_{n} + \delta_{n}$.
After applying Newton's method with this new value of the Reynolds number, we have two options:
\begin{itemize}
\item If Newton's method converges, we consider a new trial value $\operatorname{Re}_{n+2} = \operatorname{Re}_{n+1} + \delta_{n+1}$ when $\delta_{n+1} = \delta_{n}$;
\item If the method does not converge, then a halved increment $\delta_{n+1} = \frac 12 \delta_{n}$ is chosen and the calculation for the smaller trial value $\operatorname{Re}_{n+2} = \operatorname{Re}_{n} + \delta_{n+1}$ with the initial approximation to solve the problem at $\operatorname{Re}_{n}$ is performed.
\end{itemize}

The computational results for the described Reynolds number refinement are shown in Fig.~\ref{f-4}.
A total of 500 iterations are performed with $\operatorname{Re}_{1} = 500$ and $\delta_{1}=2\,000$.
The strategy of varying the incremental Reynolds number is observed most clearly on coarse grids (Figures~\ref{f-4}a,~\ref{f-4}b).
The sequential Reynolds number refinement algorithm under consideration gives us a critical value of $\operatorname{Re}$ at which an approximate solution is found.
For example, on a grid with $m=64$, we find an approximate solution at $\operatorname{Re}  \lesssim 17\,214$, on a grid with $m=128$ --- at $\operatorname{Re}  \lesssim 23\,030$, and on an even more detailed grid with $m=256$ --- at $\operatorname{Re}  \lesssim 45\,906$.
As we increase the computational grids, we find an approximate solution with larger Reynolds numbers.

\begin{figure}[htbp]
\begin{tabular}{cc}
$\operatorname{Re} = 1\,000$ & $\operatorname{Re} = 10\,000$ \\
\includegraphics[width=0.48\textwidth]{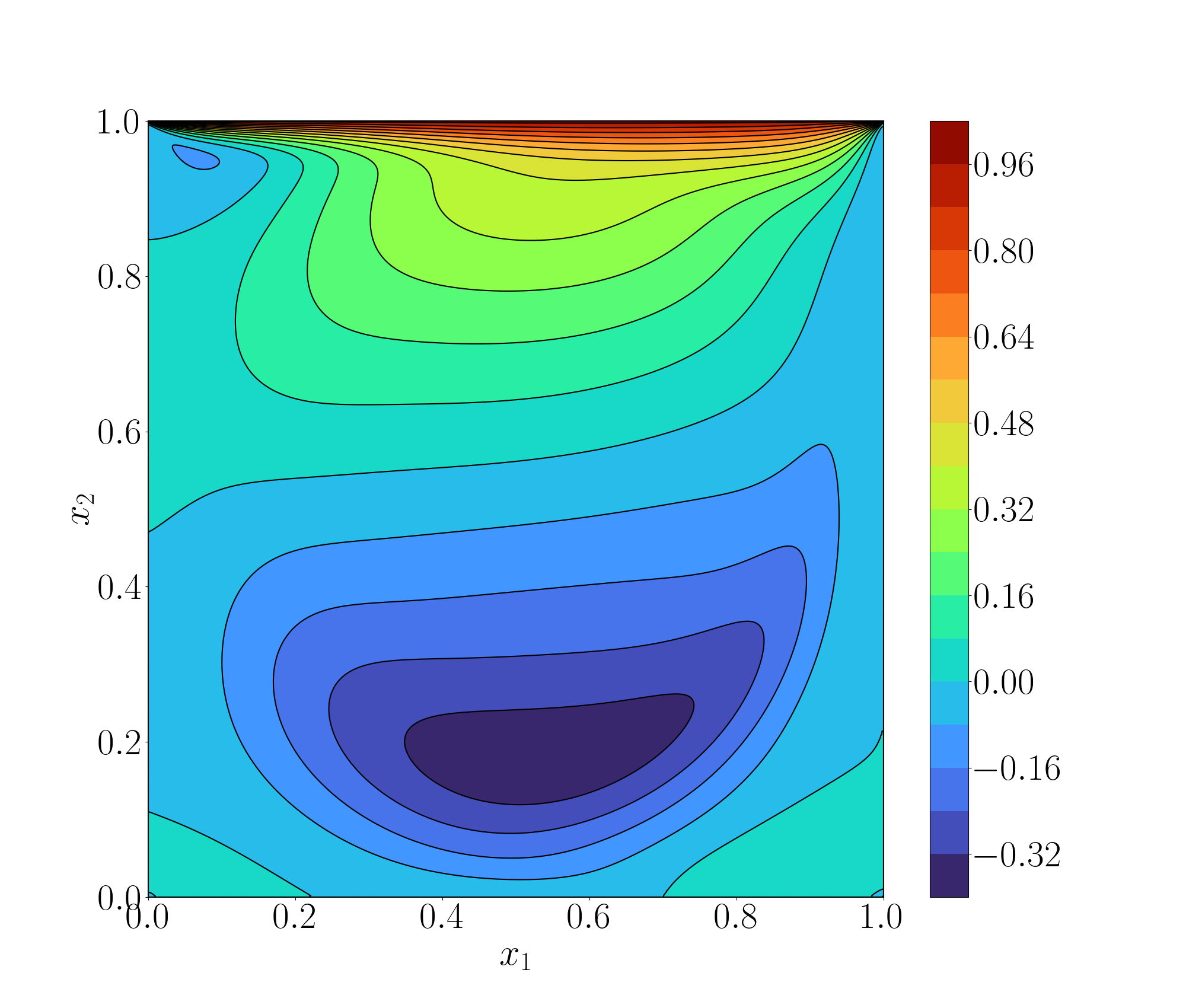}&
\includegraphics[width=0.48\textwidth]{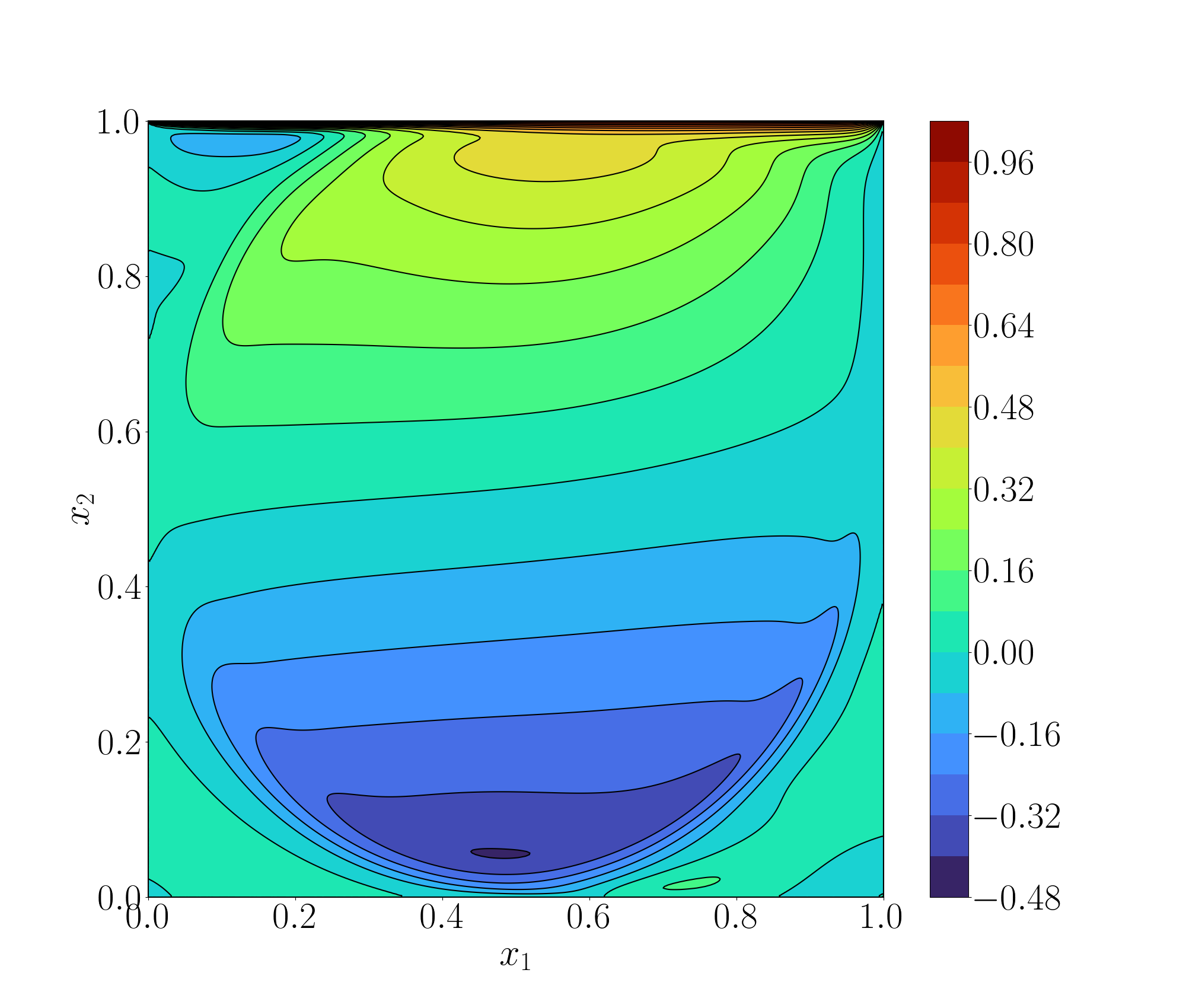}
\end{tabular}
\caption{Horizontal velocity.}
\label{f-5}
\end{figure}

\begin{figure}[htbp]
\begin{tabular}{cc}
$\operatorname{Re} = 1\,000$ & $\operatorname{Re} = 10\,000$ \\
\includegraphics[width=0.48\textwidth]{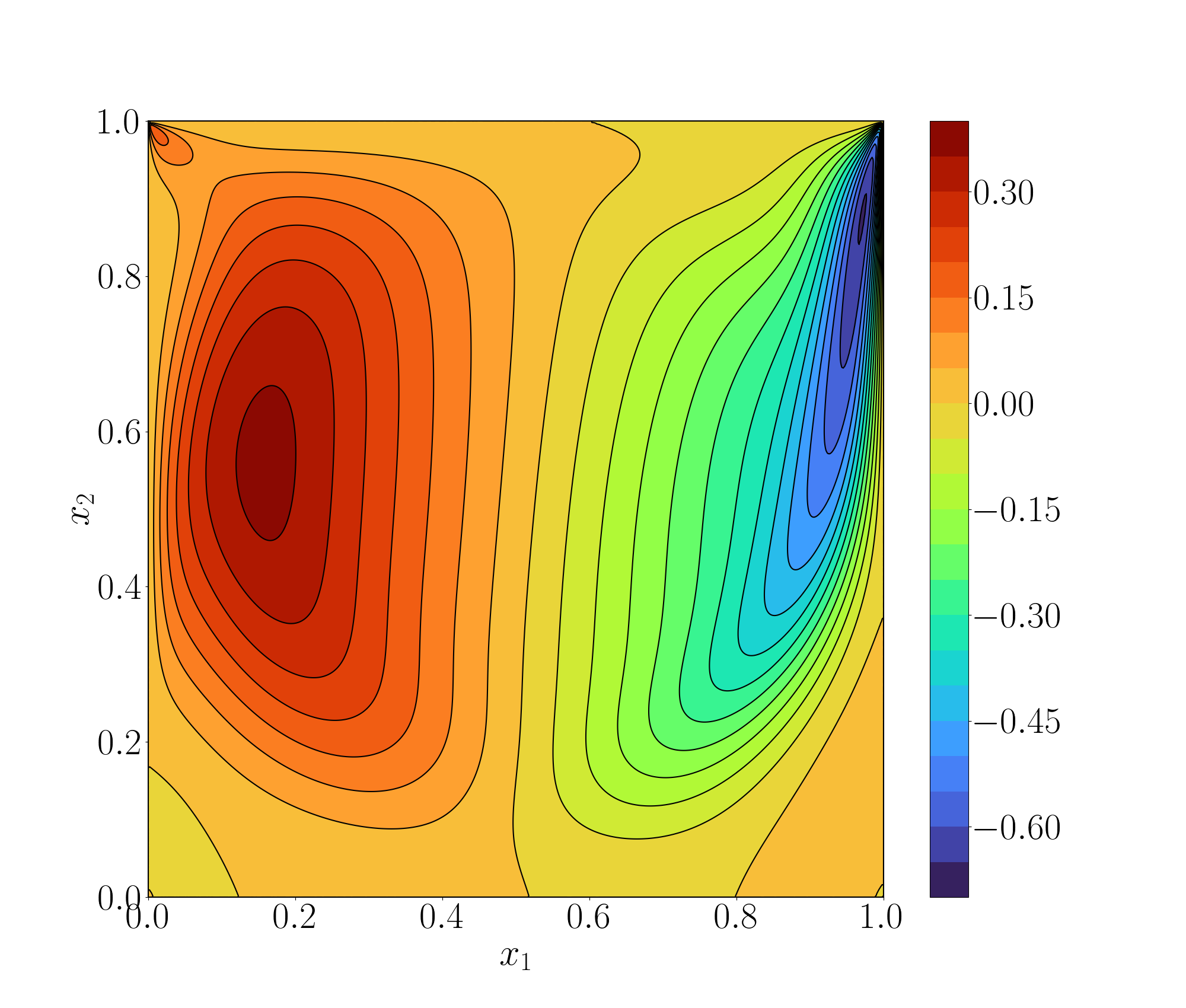}&
\includegraphics[width=0.48\textwidth]{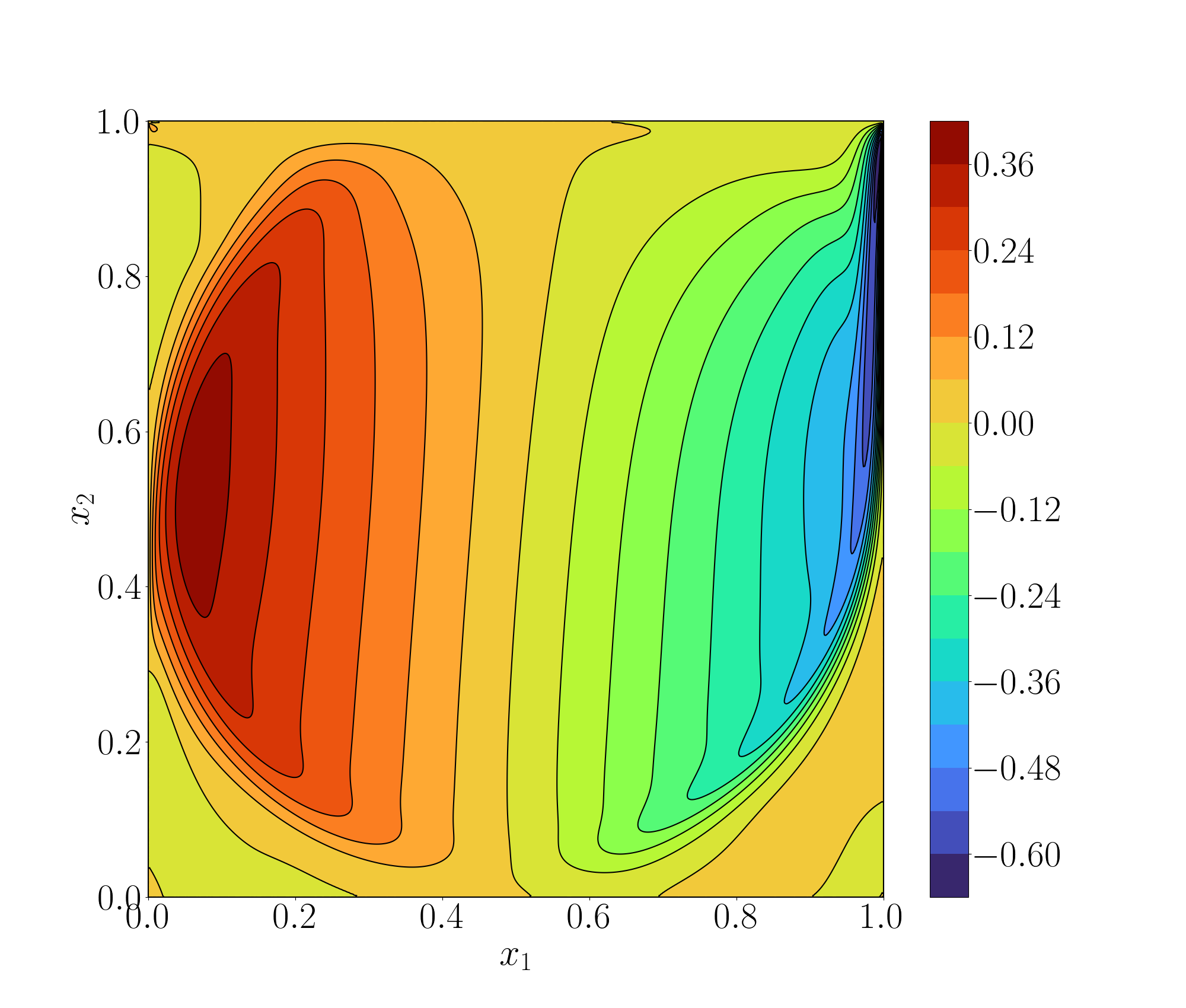}
\end{tabular}
\caption{Vertical velocity.}
\label{f-6}
\end{figure}

\begin{figure}[htbp]
\begin{tabular}{cc}
$\operatorname{Re} = 1\,000$ & $\operatorname{Re} = 10\,000$ \\
\includegraphics[width=0.48\textwidth]{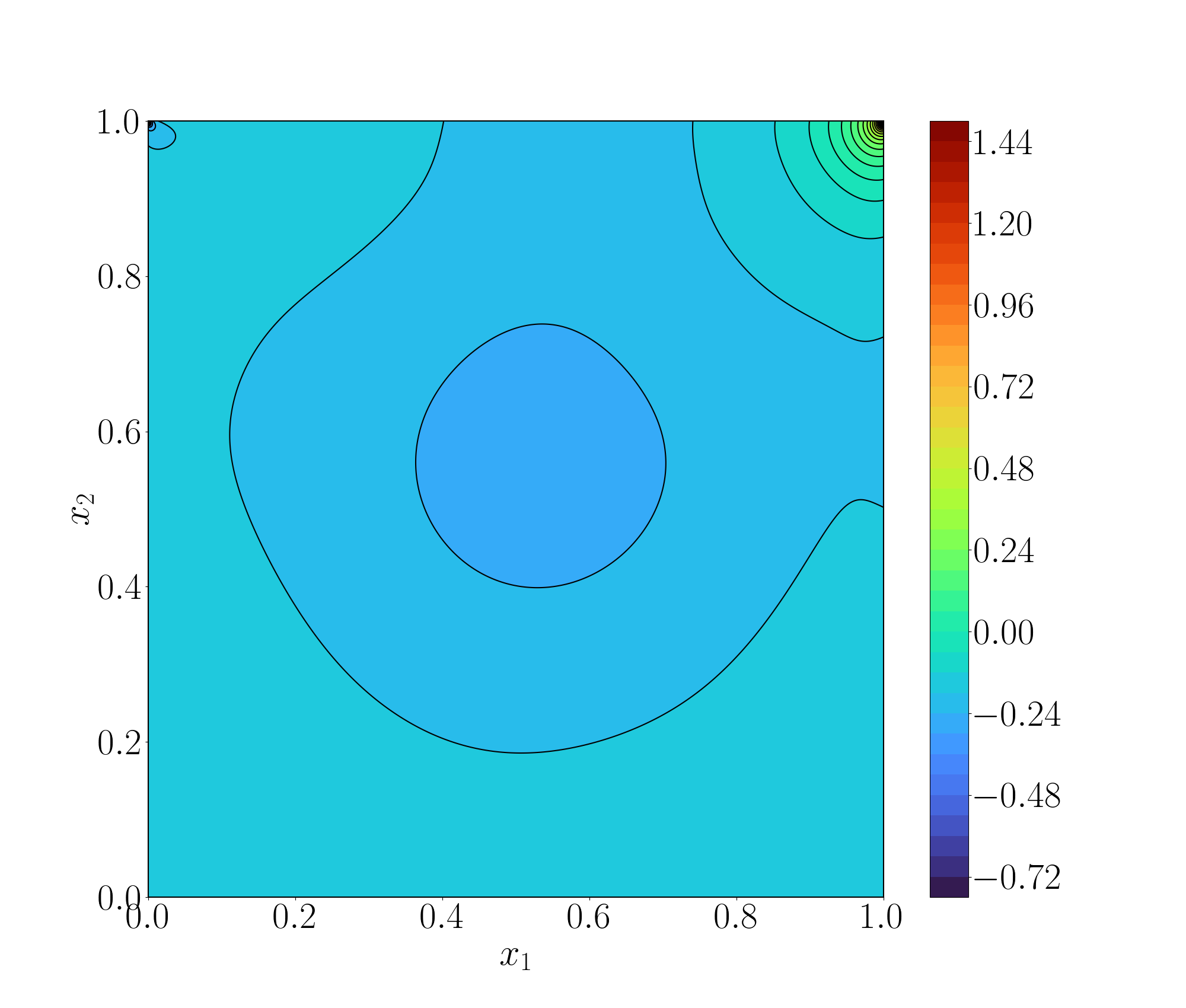}&
\includegraphics[width=0.48\textwidth]{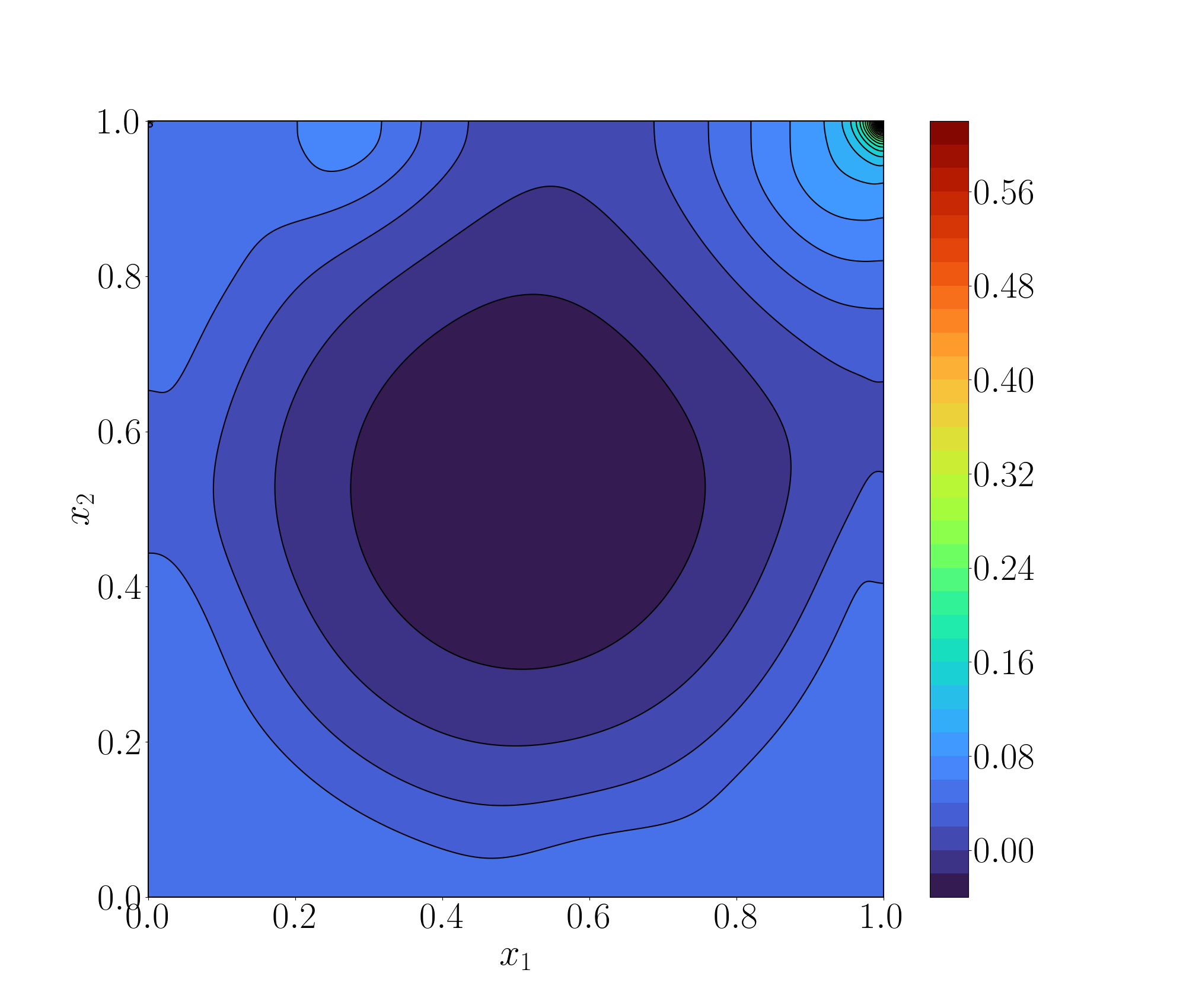}
\end{tabular}
\caption{Pressure.}
\label{f-7}
\end{figure}

\begin{figure}[htbp]
\begin{tabular}{cc}
$\operatorname{Re} = 1\,000$ & $\operatorname{Re} = 10\,000$ \\
\includegraphics[width=0.48\textwidth]{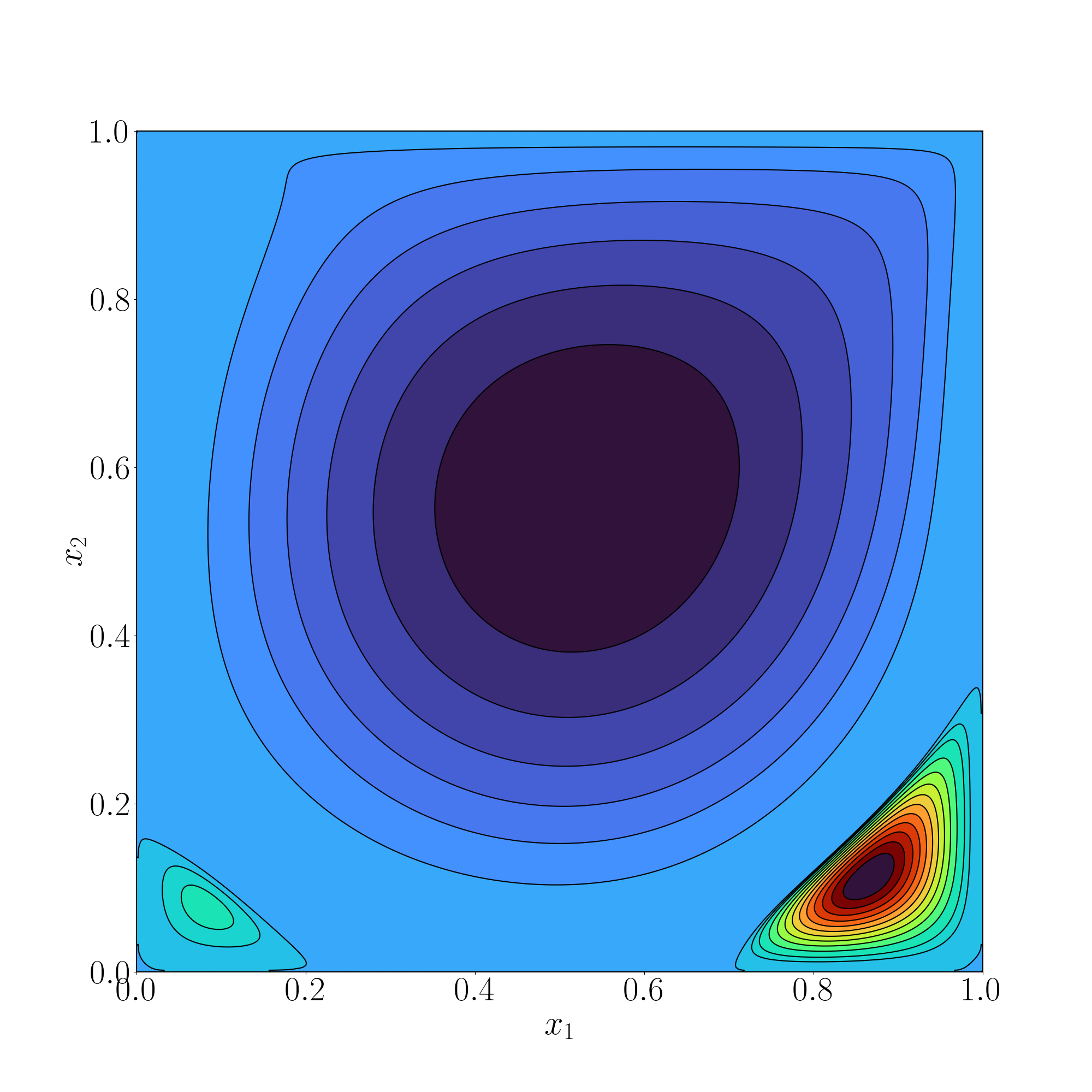}&
\includegraphics[width=0.48\textwidth]{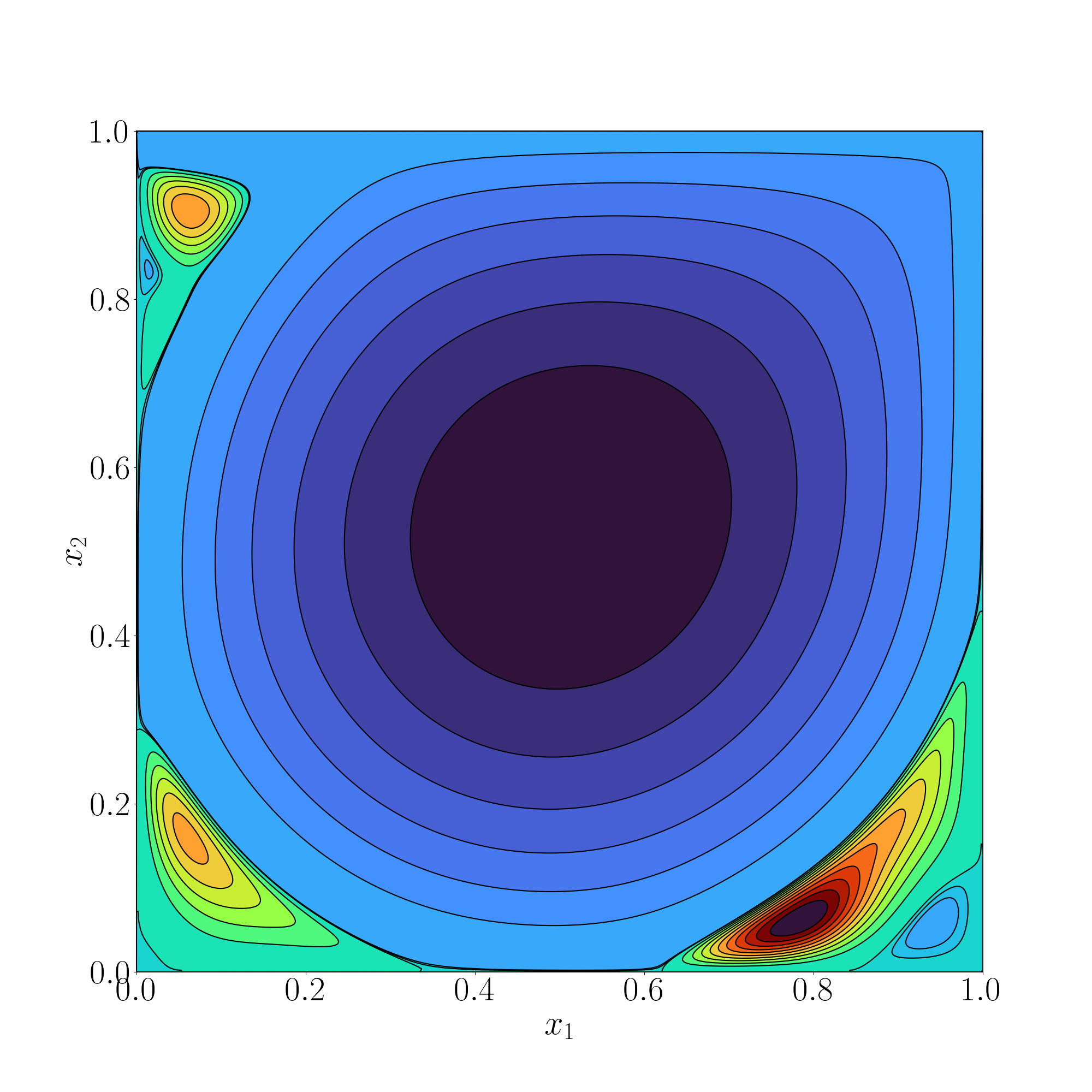}
\end{tabular}
\caption{Stream function.}
\label{f-8}
\end{figure}

\begin{figure}[htbp]
\begin{tabular}{cc}
$\operatorname{Re} = 1\,000$ & $\operatorname{Re} = 10\,000$ \\
\includegraphics[width=0.48\textwidth]{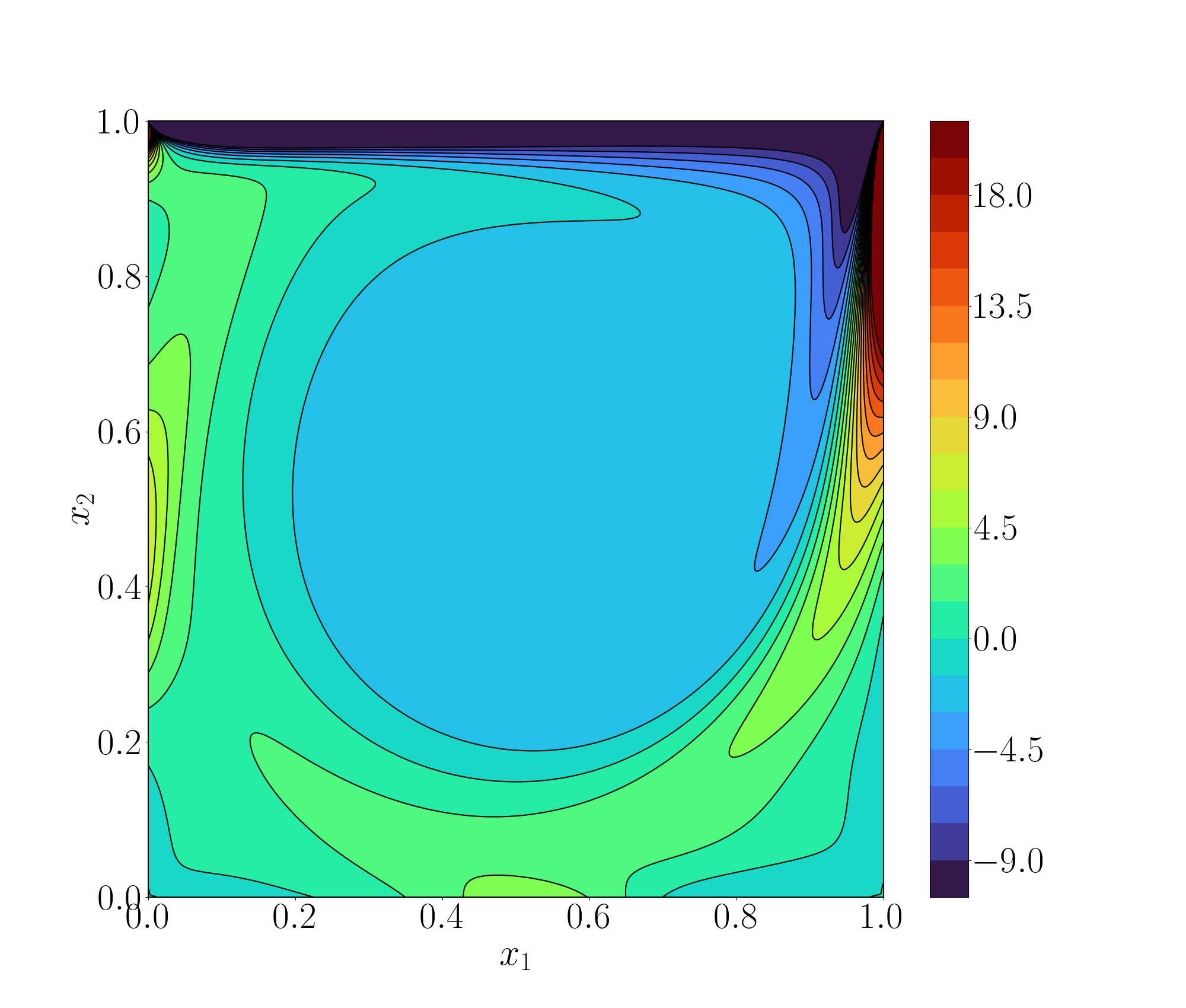}&
\includegraphics[width=0.48\textwidth]{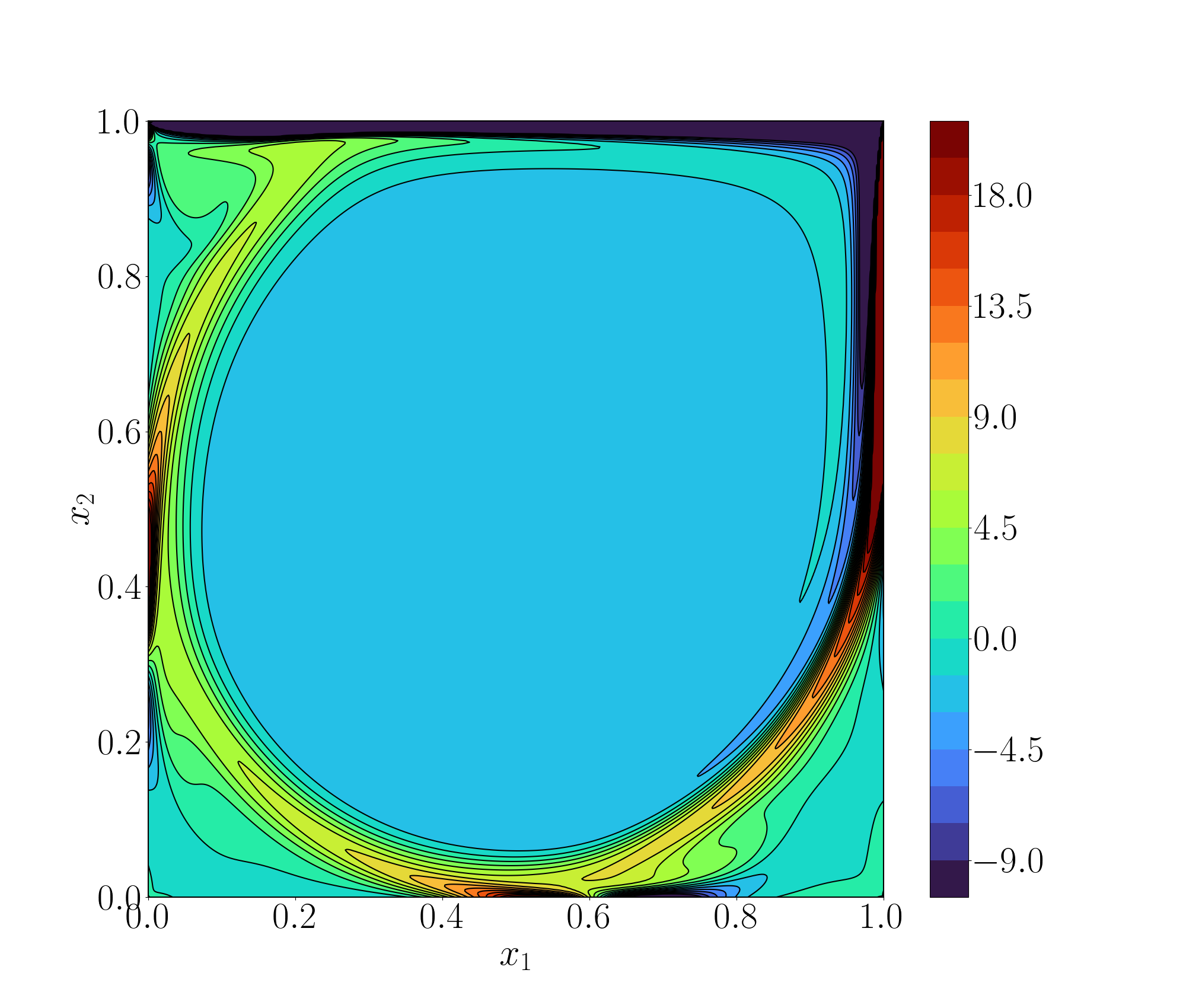}
\end{tabular}
\caption{Vorticity.}
\label{f-9}
\end{figure}

The approximate solution at Reynolds numbers of $1\,000$ and $10\,000$ is shown in Fig.~\ref{f-5}-\ref{f-9}.
A uniform triangular grid with $m=512$ was used.

\clearpage

\section{Other iterative methods}

Let us consider the possibility of constructing other iterative methods besides Newton's method, for which a good initial approximation is necessary for convergence.
We would like the iterative method to be less dependent on the initial condition than Newton's method while conceding the speed of convergence.

\subsection{Iterative methods with relaxation}

We will stay in the class of two-level iterative methods, where the new approximation $\{\bm u^{k+1}, p^{k+1}\}$ is computed from the solution at the previous iteration $\{\bm u^{k}, p^{k}\}$.
The choice of the iterative method is related to the linearization of the convective term in Eq. (\ref{1}).
Given the bilinearity of $C (\bm v, \bm u)$, we will associate the linearizations with the choice of $\bm v$ or $\bm u$ in the previous iteration.

It is natural to associate such linearizations with the relaxation of the approximate solution at a new iteration.
In this case, similar to (\ref{5})--(\ref{7}) one finds an auxiliary solution $\{\widetilde{\bm u}^{k+1}, \widetilde{p}^{k+1}\}$ of the boundary value problem
\begin{equation}\label{11}
\begin{split}
\big (C ( \widetilde{\bm u}\otimes \widetilde{\bm u}) \big )^{k+1} + \operatorname{grad} \widetilde{p}^{k+1}
& - \frac{1}{\operatorname{Re}} \operatorname{div} \operatorname{grad}  \widetilde{\bm u}^{k+1} = 0 , \\
\operatorname{\operatorname{div}}  \widetilde{\bm u}^{k+1} & = 0,
\quad \bm x \in \Omega , \\
\widetilde{\bm u}^{k+1}(\bm x) & = \bm g(\bm x),
\quad \bm x \in \partial \Omega .
\end{split}
\end{equation}
The solution at the new iteration is taken as
\[
\begin{split}
\bm u^{k+1} = \sigma \bm {\widetilde{u}}^{k+1} + (1-\sigma)\bm u^{k} , \\
p^{k+1} = \sigma {\widetilde{p}}^{k+1} + (1-\sigma) p^{k} ,
\end{split}
\]
where $\sigma$ is relaxation parameter.

Method 1 is based on linearization by the first argument of the convective summand: $C (\bm v^k, \bm u^{k+1})$.
In this case, in (\ref{11}) we have
\begin{equation}\label{12}
\big (C ( \widetilde{\bm u}\otimes \widetilde{\bm u}) \big )^{k+1} = C ( \widetilde{\bm u}^k\otimes \widetilde{\bm u}^{k+1})
= \operatorname{div}(\widetilde{\bm u}^{k} \otimes \widetilde{\bm u}^{k+1} ).
\end{equation}
Method 2 performs linearization on the second argument: $C (\bm v^k, \bm u^{k+1})$.
In the first equation (\ref{11}), we assume
\begin{equation}\label{13}
\big (C ( \widetilde{\bm u}\otimes \widetilde{\bm u}) \big )^{k+1} = C ( \widetilde{\bm u}^{k+1} \otimes \widetilde{\bm u}^{k})
= \operatorname{div}(\widetilde{\bm u}^{k+1} \otimes \widetilde{\bm u}^{k} ).
\end{equation}
Method 3 combines the linearizations of method 1 and method 2 when
\begin{equation}\label{14}
\begin{split}
\big (C ( \widetilde{\bm u}\otimes \widetilde{\bm u}) \big )^{k+1} & = \frac 12 \Big (C ( \widetilde{\bm u}^k\otimes \widetilde{\bm u}^{k+1})
+ C ( \widetilde{\bm u}^{k+1} \otimes \widetilde{\bm u}^{k}) \Big ) \\
& = \frac 12 \big (\operatorname{div}(\widetilde{\bm u}^{k} \otimes \widetilde{\bm u}^{k+1} ) 
+ \operatorname{div}(\widetilde{\bm u}^{k+1} \otimes \widetilde{\bm u}^{k} ) \big ).
\end{split}
\end{equation}
Method 4 is nothing but Newton's method.
In this case, linearization is performed by the rule (\ref{8}).

\subsection{Numerical experiments}

\begin{figure}[htbp]
\begin{tabular}{cc}
Method 1 & Method 2 \\
\includegraphics[width=0.45\textwidth]{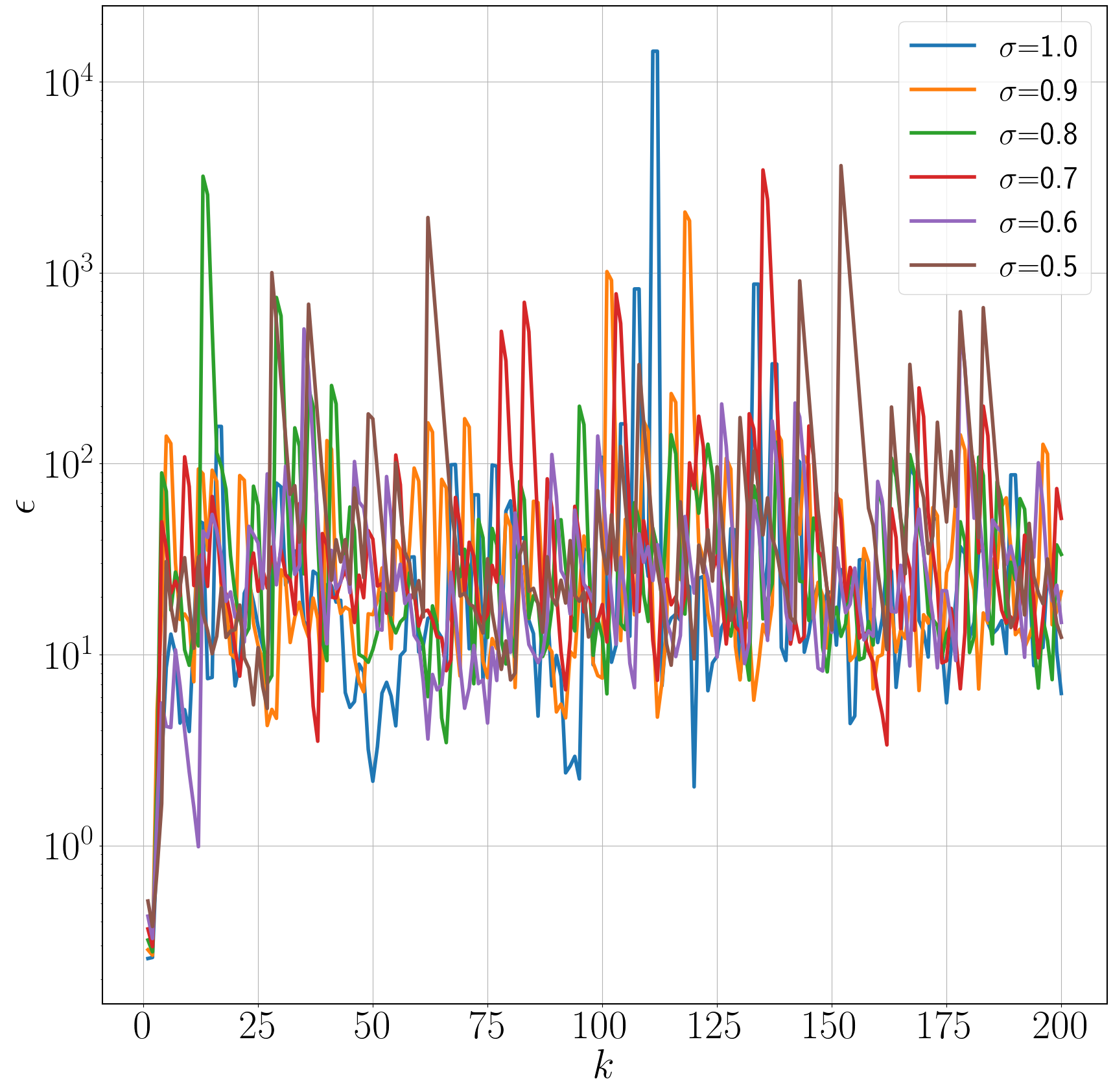} &
\includegraphics[width=0.45\textwidth]{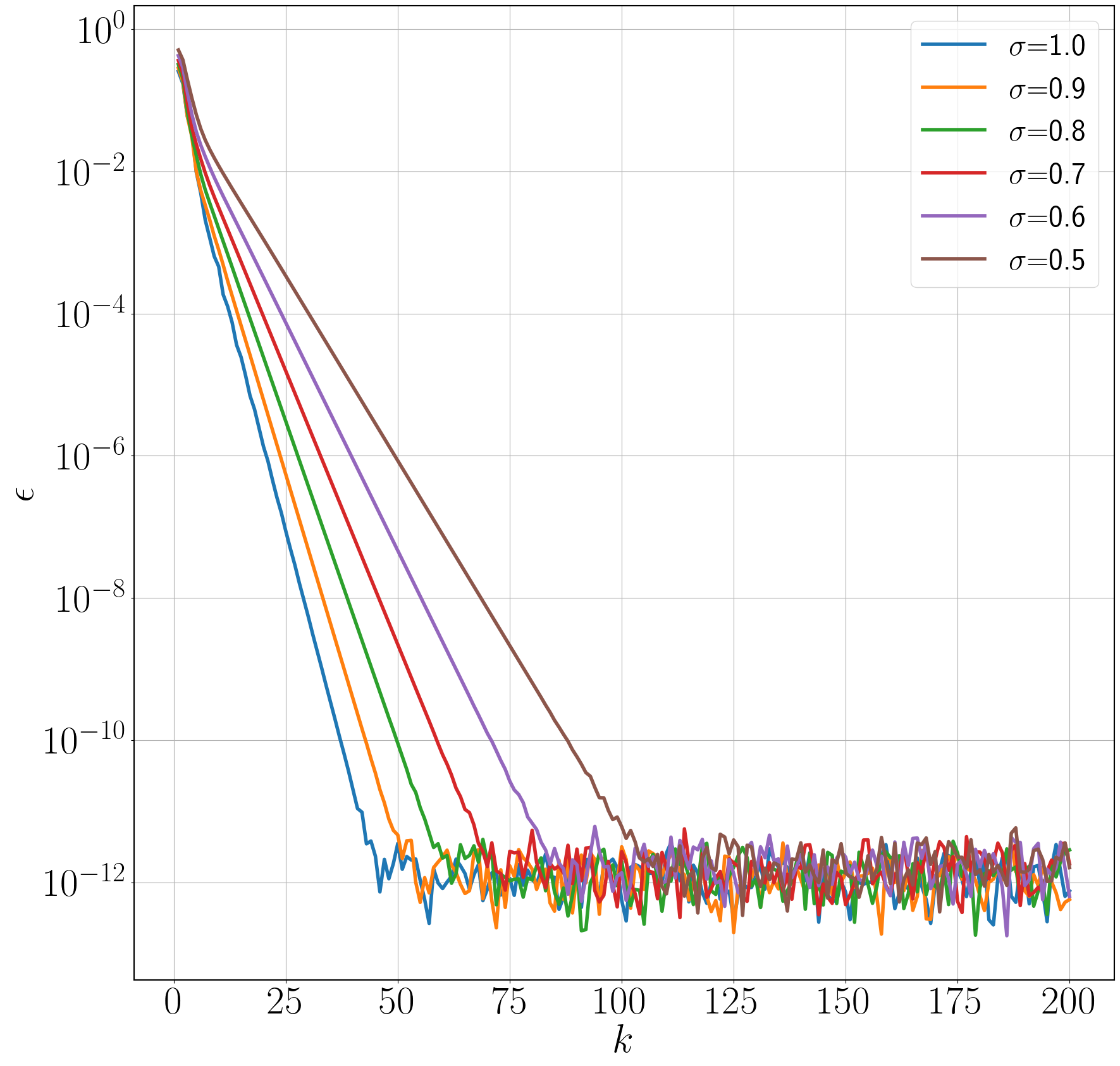}\\
Method 3 &  Method 4 \\
\includegraphics[width=0.45\textwidth]{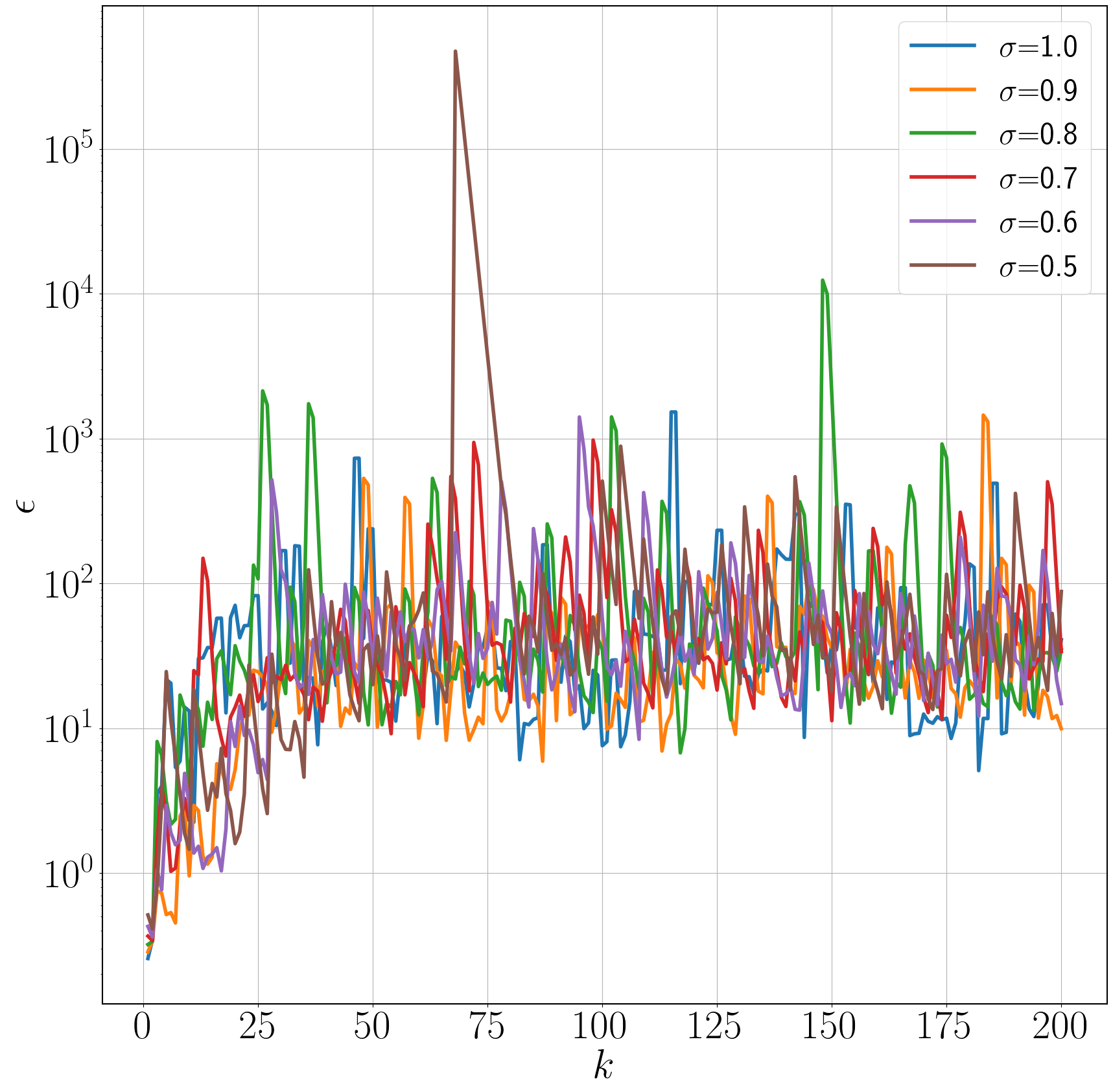} &
\includegraphics[width=0.45\textwidth]{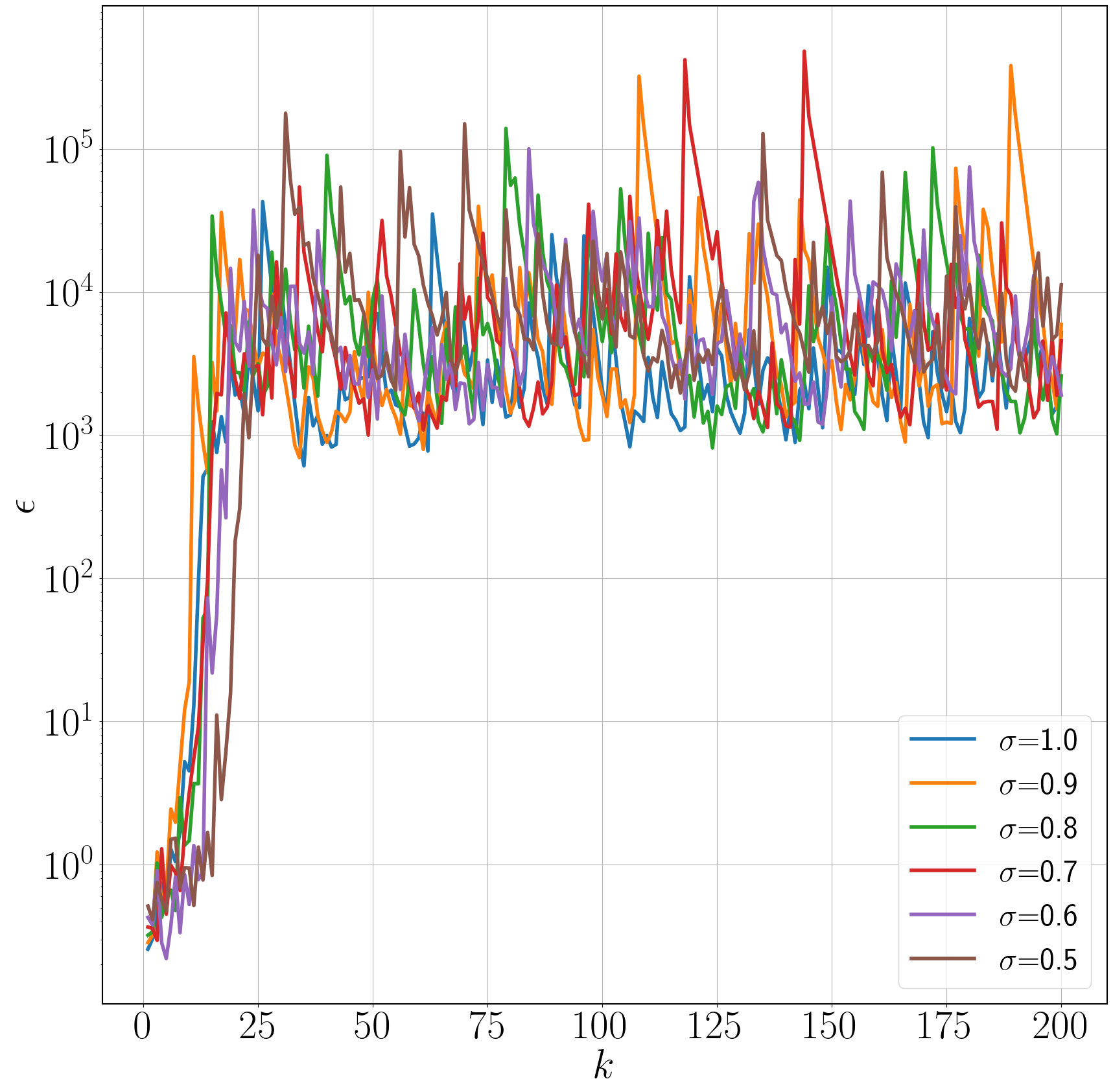}
\end{tabular}
\caption{Convergence of the iterative method with relaxation under different linearizations of the convective term.}
\label{f-10}
\end{figure}

Calculations are performed on a grid with $m = 128$ at $\operatorname{Re} = 1\,000$.
The initial approximation is taken as $\bm u^0 (\bm x) = 0, \ \bm x \in \Omega$.
The computational results using the iterative method (\ref{11}) with linearizations (\ref{12})--(\ref{14}) as well as (\ref{8}) are shown in Fig.~\ref{f-10}.
The main conclusion is that only iterative method 2 converges at $\sigma \leq 1$.
Other methods do not converge, and reducing the relaxation parameter does not help.

\begin{figure}[htbp]
\begin{tabular}{cc}
$\operatorname{Re} = 5\,000$ & $\operatorname{Re} = 10\,000$ \\
\includegraphics[width=0.48\textwidth]{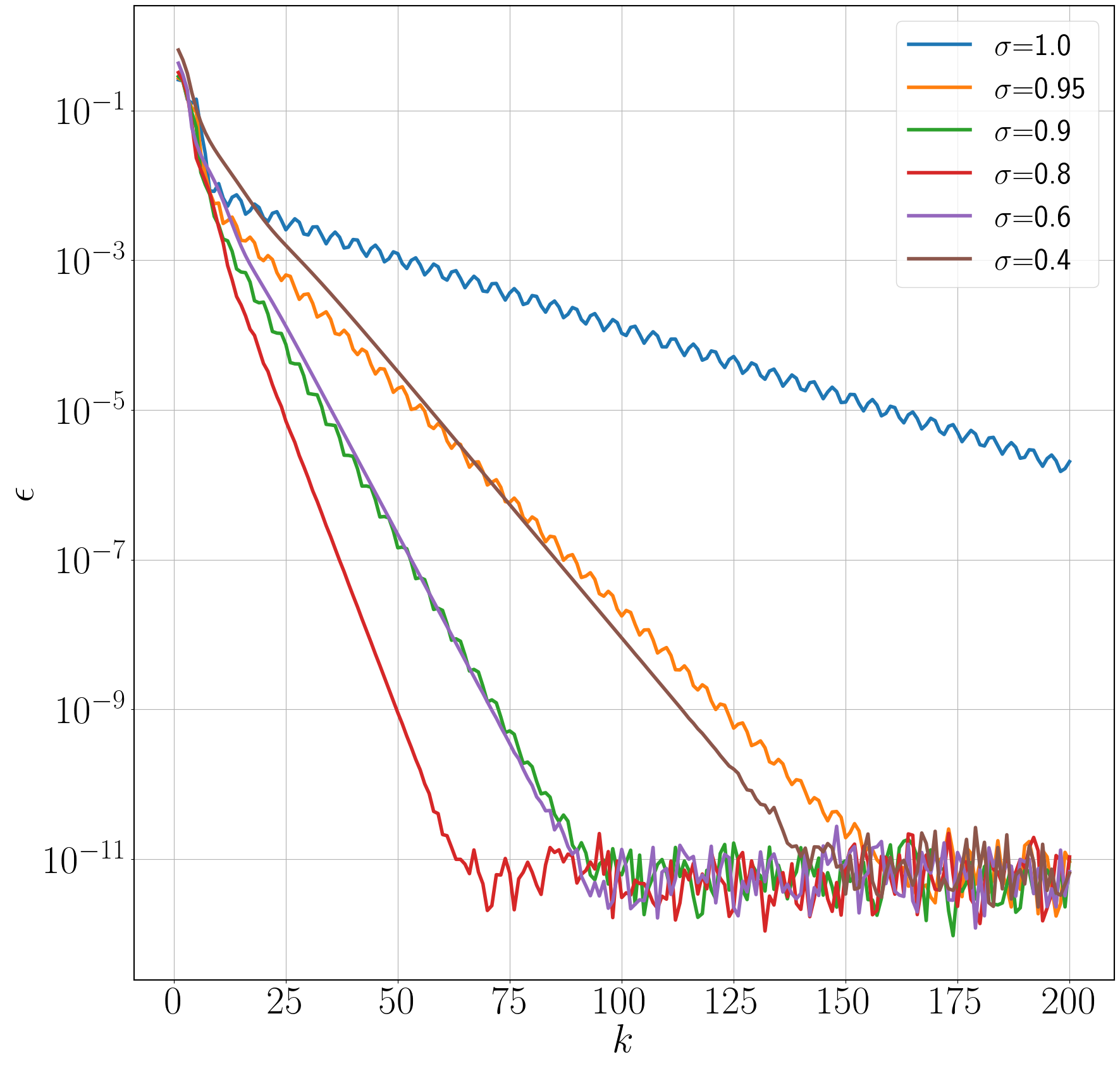}&
\includegraphics[width=0.48\textwidth]{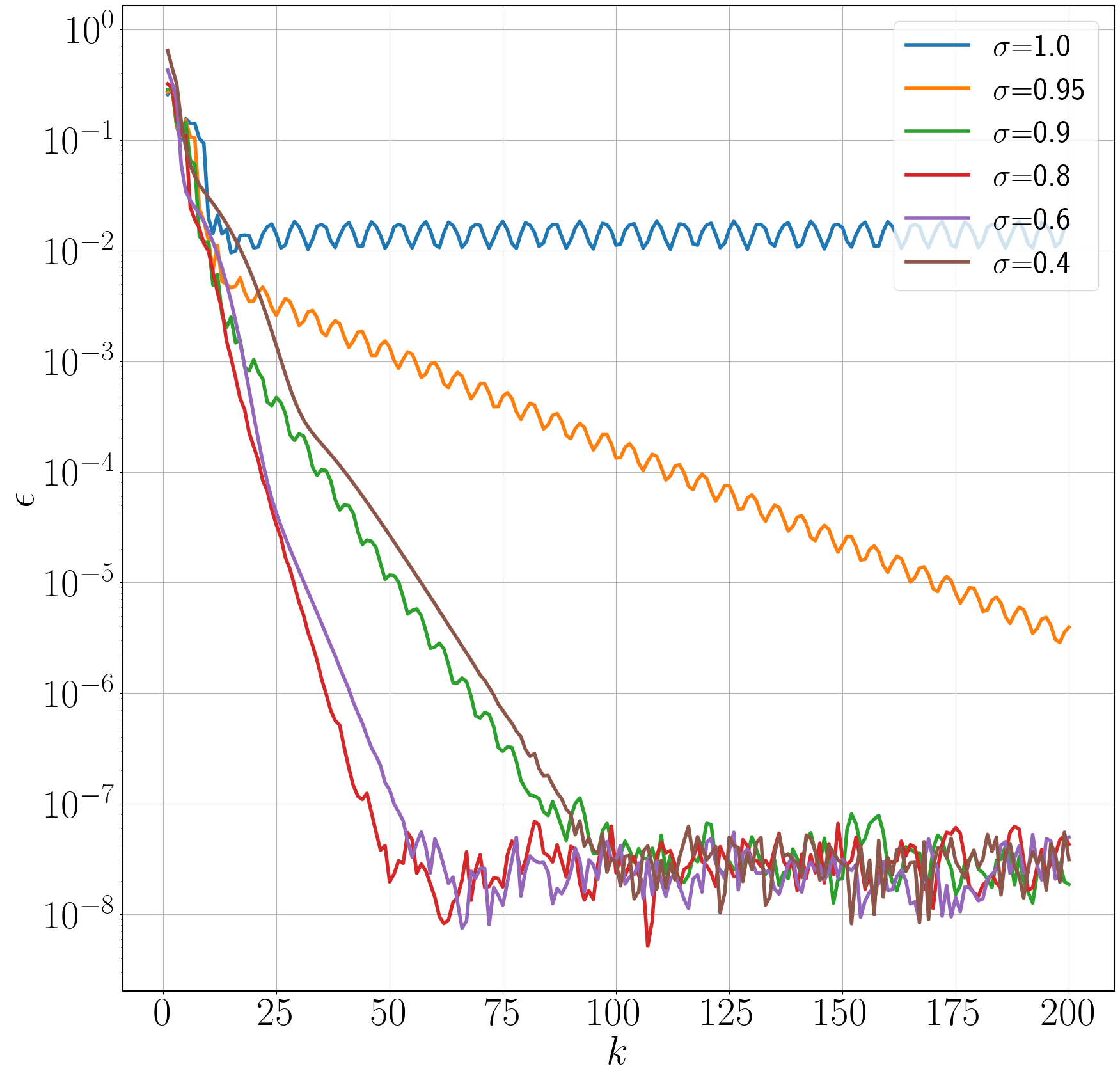}\\
\end{tabular}
\caption{Convergence of method two at different Reynolds numbers.}
\label{f-11}
\end{figure}

Method 2 converges linearly, while Newton's method if it converges, converges quadratically.
Method 2 has the advantage of a much larger region of convergence when given an initial approximation.
The possibilities of using this method for approximate solutions to problems with large Reynolds numbers are illustrated in Fig.~\ref{f-11}.
We consider the flow problem in a square cavity at $\operatorname{Re} = 5\,000$ and $\operatorname{Re} = 10\,000$.
Convergence is observed at $\sigma \leq \sigma_0(\operatorname{Re})$, where $\sigma_0(\operatorname{Re})$ decreases as $\operatorname{Re}$ increases.

\subsection{Multiple solutions}

For large values of the Reynolds number, we cannot have one solution but several.
Computational algorithms for investigating such a bifurcation problem for solutions to nonlinear problems are currently attracting much attention from researchers.
Successes in obtaining multiple solutions for stationary problems of incompressible fluid hydrodynamics are reported in a recent paper \cite{erturk2022bifurcation}.
Let us note the capabilities of the iterative method 2 discussed above.

As a test problem, following the paper \cite{erturk2022bifurcation}, we will consider
the boundary value problem for the Navier-Stokes equations for an incompressible fluid inside a semi-elliptical lid-driven cavity and an ellipse axis ratio of 2:1.
The solution is non-unique for Reynolds numbers above the critical Reynolds number equal to about $5\,059$.

\begin{figure}[htbp]
\begin{tabular}{c}
$\sigma=0.5$ \\
\includegraphics[width=1\textwidth]{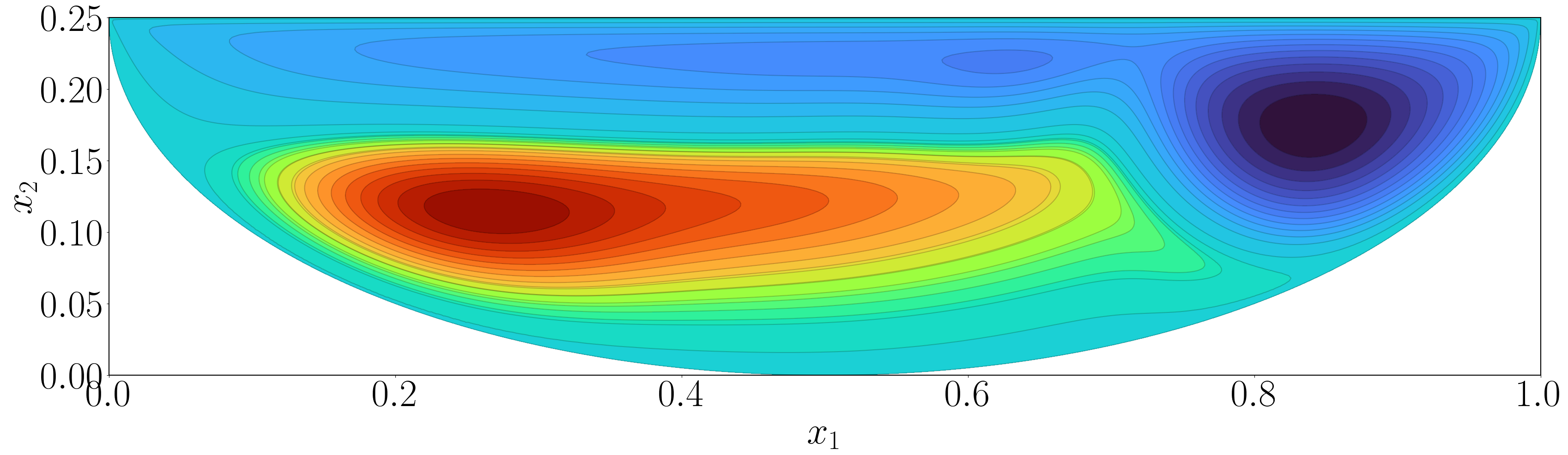} \\
$\sigma=0.7$ \\
\includegraphics[width=1\textwidth]{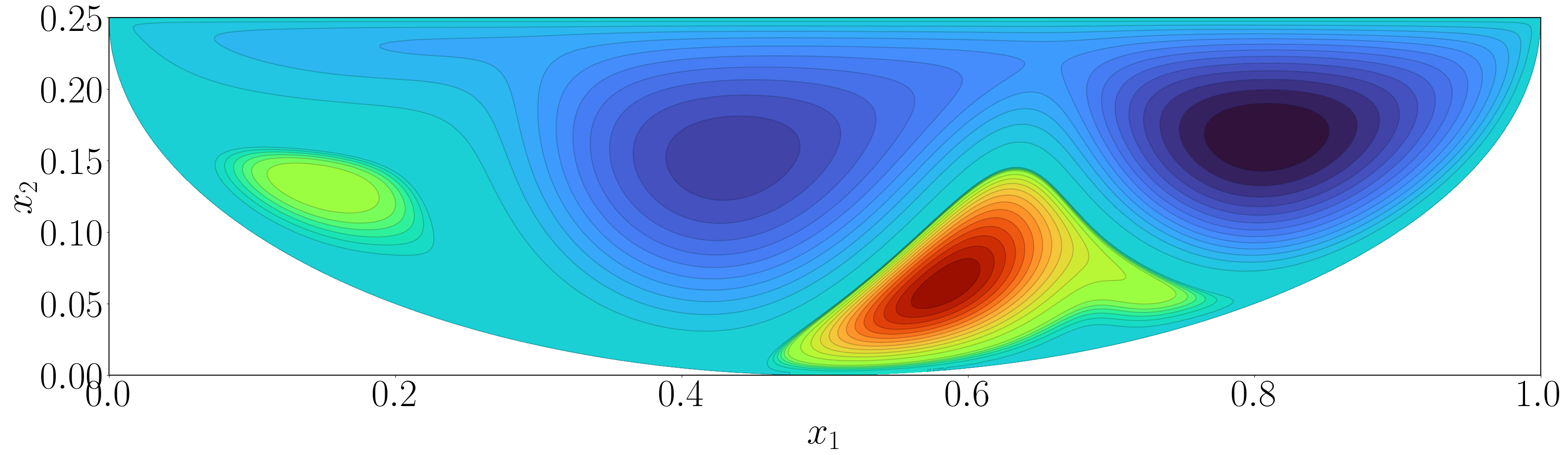}\\
\end{tabular}
\caption{Stream function, $\operatorname{Re}=5\,500$.}
\label{f-12}
\end{figure}

\begin{figure}[htbp]
\begin{tabular}{c}
$\sigma=0.5$ \\
\includegraphics[width=1\textwidth]{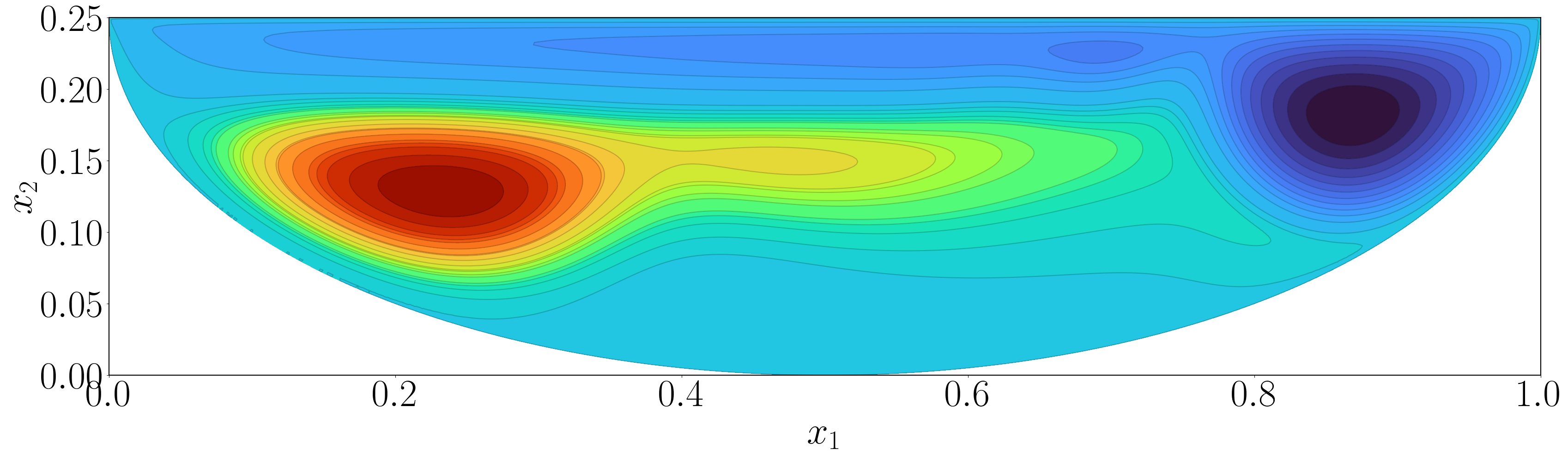} \\
$\sigma=0.7$ \\
\includegraphics[width=1\textwidth]{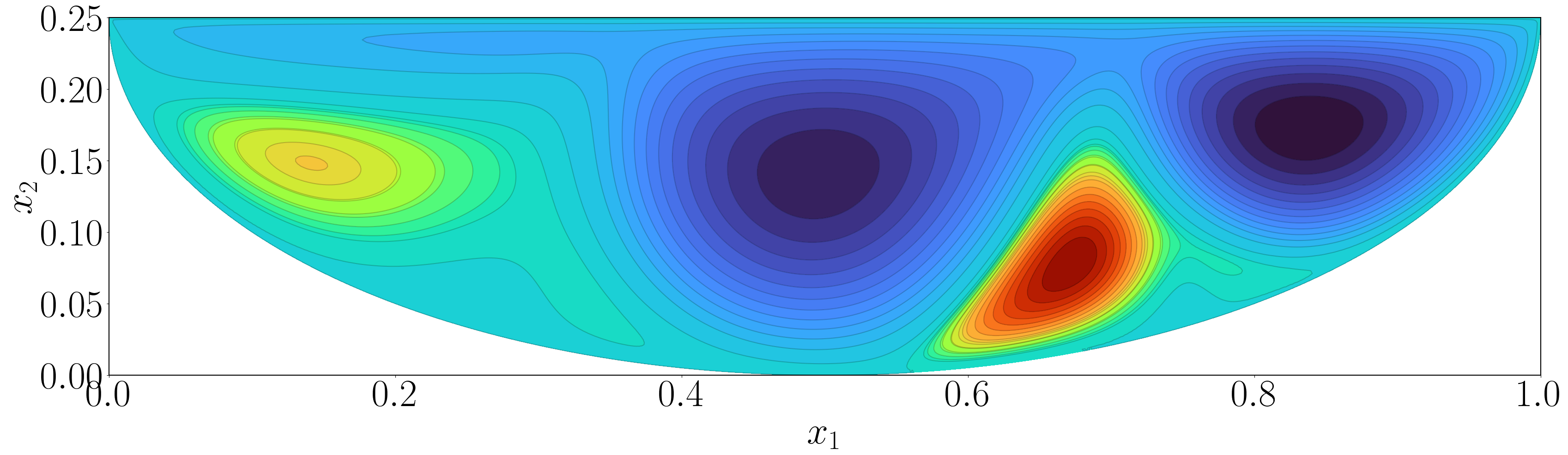}\\
\end{tabular}
\caption{Stream function, $\operatorname{Re}=10\,000$.}
\label{f-13}
\end{figure}

We use a triangular computational grid with $ 228\,348$ nodes and an initial condition $\bm u^0 (\bm x) = 0, \ \bm x \in \Omega$.
When the critical Reynolds number is exceeded, we obtain either one or the other approximate solution.
Two qualitatively different solutions at $\operatorname{Re}=5\,500$ are shown in Fig.~\ref{f-12}.
When the relaxation parameter $\sigma = 0.7$, we obtain the first solution with one large negative vortex in the right part and a large positive vortex in the left. At $\sigma = 0.9$, we obtain the second solution, characterized by two large negative vortexes and a large positive vortex between them.

Similar data at higher Reynolds numbers are shown in Fig.~\ref{f-5}. We observe a change in the size of negative vortices with increasing $\operatorname{Re}$ for both solutions. The left negative vortex increases in size in the second solution, and the right one decreases for both solutions. The same value of $\sigma$ can give both the first and second solutions. When $\operatorname{Re} =5\,500$ and $\sigma=0.7$ we get the first solution, when $\operatorname{Re} =10\,000$ and $\sigma=0.7$ --- the second solution.

\begin{figure}[htbp]
\begin{tabular}{cc}
$\operatorname{Re} = 5\,500$ & $\operatorname{Re} = 10\,000$ \\
\includegraphics[width=0.48\textwidth]{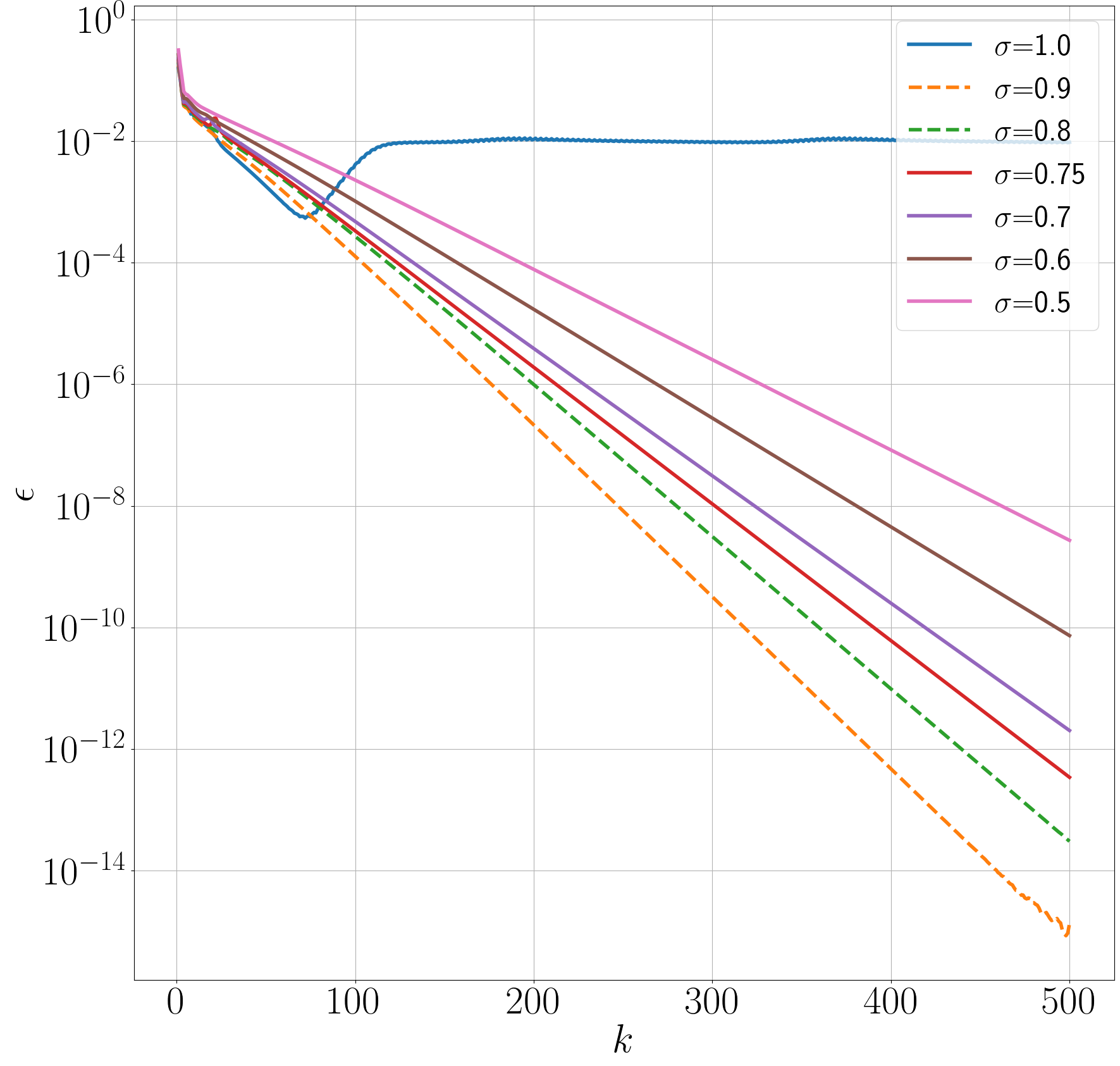}&
\includegraphics[width=0.48\textwidth]{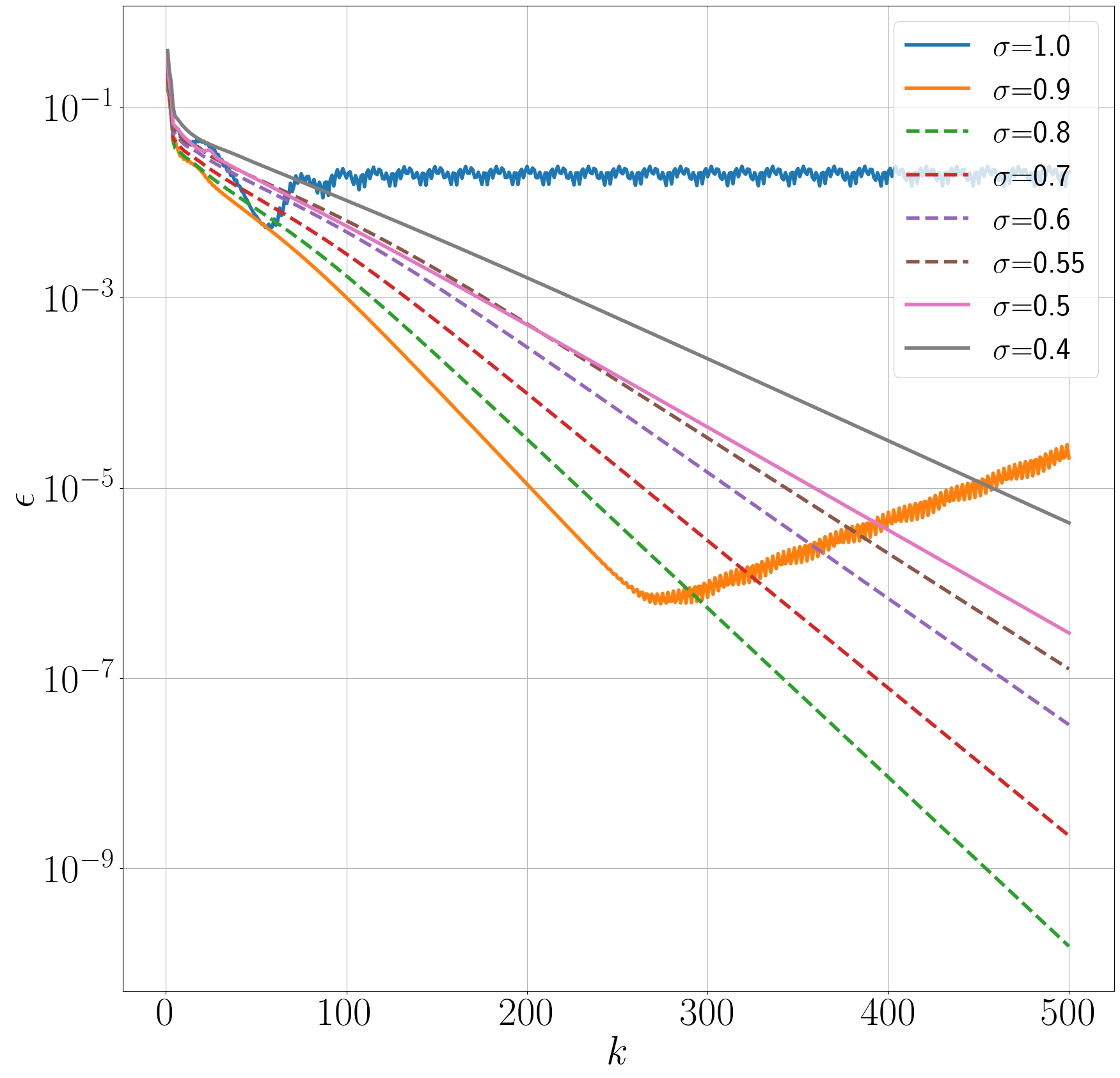}\\
\end{tabular}
\caption{Convergence to different solutions for different values of the relaxation parameter.}
\label{f-14}
\end{figure}

The influence of the relaxation parameter is traced in more detail in Fig.~\ref{f-14}.
For the problem with $\operatorname{Re} = 5\,500$, iterative method 2 fails to converge at $\sigma = 1$.
For $\sigma = 0.9, 0.8$ we have the first solution, and for $\sigma = 0.75, 0.7, 0.6, 0.5$ we get the second solution (see Fig.~\ref{f-12}).
When $\operatorname{Re} = 10\,000$, iterative method 2 fails to converge at $\sigma = 1, 0.9$.
We have the first solution (see Fig.~\ref{f-13}) at $\sigma = 0.8, 0.7, 0.6, 0.55$, and the second solution --- at $\sigma = 0.5, 0.4$.

\section{Conclusions}

\begin{enumerate}[(1)]
\item For approximate solution of stationary boundary value problems for the Navier-Stokes equations for incompressible fluid in pressure-velocity variables, iterative methods related to different linearizations of the convective component are proposed.
Newton's method at zero initial approximation allows for the obtaining of an approximate solution only at small Reynolds numbers.
In general, the solution is found by sequentially increasing the Reynolds number.
\item Numerical experiments are carried out for two-dimensional fluid flows in a cavity using different linearizations of the nonlinear convective summand in the Navier-Stokes equations. We have constructed an iterative method with relaxation that gives a solution at large Reynolds numbers with zero initial approximation.
\item Using the proposed iterative method, one solution for Reynolds numbers below the critical Reynolds number and two for Reynolds numbers above the critical Reynolds number is obtained for flows in a cavity of the semi-elliptical cross-section. The two solutions are obtained by choosing different values of the relaxation parameter.
\end{enumerate} 


\end{document}